\documentclass[12pt]{article}
\usepackage{graphicx}

\newcommand{\C}{{\bf C}}

\newcommand{\R}{{\bf R}}
\newcommand{\T}{{\bf T}}

\newcommand{\Z}{{\bf Z}}

\def\gato{\sharp}
\def\sc{\mathop{\gato}\limits}
\def\scf{\mathop{\amalg}\limits}

\begin{document}

\begin{center}

\textbf{INTERSECTIONS OF QUADRICS, MOMENT-ANGLE MANIFOLDS AND CONNECTED SUMS.}

\bigskip

\textbf{Samuel Gitler.}

\smallskip

\textbf{Santiago L\'opez de Medrano.}

\end{center}

\textit{
\begin{flushright}
To Sylvain Cappell, on his 65th birthday.
\end{flushright}
}
\bigskip

\section*{Introduction.}

\bigskip

The geometric topology of the generic intersection of two homogeneous coaxial quadrics in $\R^n$ was studied in [LdM1] where it was shown that its intersection with the unit sphere is in most cases diffeomorphic to a triple product of spheres or to the connected sum of sphere products. The proof involved a geometric description of the group actions on them and of their polytope quotients as well as a splitting of the homology groups of those manifolds. It also relied heavily on a normal form for them and many related computations. The part about group actions, polytopes and homology splitting was equally valid for the intersection of any number of such quadrics, but
the obstacle to extending the main result for more than two turned out to be the hopeless-looking problem of finding their normal forms, close to that of classifying all simple polytopes.

\bigskip

The study of those intersections continued in other directions, especially those related to the projectivizations of their complex versions (known now as \emph{LV-M manifolds}, see [LdM2], [LdM-V], [Me], [L-N], [B], [Me-V1], [B-M] and the recent review [Me-V2]), giving many new examples of non-algebraic complex manifolds which fibre over toric varieties.

\bigskip

Following this line of research, in [B-M] a deep study of LV-M manifolds included also important advances in the study of the geometric topology of the intersection of $k$ homogeneous quadrics in their complex versions (the manifolds now called \emph{moment-angle manifolds}). The main questions addressed were the following:

\bigskip

1) Whether they can always be built up from spheres by repeatedly taking products or connected sums: they produced new examples for any $k$ which are so, but also showed how to construct many cases which are not. Many interesting questions arose, including a specific conjecture.

\bigskip

2) The transition between different topological types when the generic condition is broken at some point of a deformation (\textit{wall-crossing}).

\bigskip

3) A product rule of their cohomology ring (in the spirit of the description of the homology of $Z$ given in [LdM1]) and its applications to question 1).

\bigskip

Meanwhile, and independently, in [D-J] essentially the same manifolds were constructed in a more abstract way, where the main objective was to study the algebraic topology of some important quotients of them called initially \emph{toric manifolds} and now \emph{quasitoric manifolds}. This article originated an important development through the work of many authors, and there is a vast and deep literature along those lines for which the reader is referred to [B-P]. Yet for a long time no interchange occurred between the two lines of research involving the same objects, until small connections appeared in the final version of [B-M] and in [D-S]. In particular, it turned out that examples relevant to question 1) above were known to these authors (see [B2]), and in [B1] Baskakov had a product rule for the cohomology ring, similar but dual to that of [B-M] mentioned in 3) above. All those examples were independent and more or less simultaneous, yet both product rules followed from an earlier computation by Buchstaber and Panov of the cohomology ring\footnote{In the first version of [B-M] the product rule was derived from results in [dL], while in the published version it relied on [B-P].}.

\bigskip

One recent expression of the line of research derived from [D-J] is the article [B-B-C-G] where a far-reaching generalization is made and a general geometric splitting formula is derived that is, in particular, a very good tool for understanding the relations among the homology groups of different spaces. This understanding turned out to be fundamental for us in tracing a way
through the abstract situation of the intersection of $k>2$ quadrics, thus combining efficiently both approaches to the subject as expressed in [LdM1] and [B-B-C-G] to obtain the results in the present article. Nevertheless, these results do not depend logically on those articles and are actually more geometric than any of them, involving practically no computations.

\bigskip

The results in this article follow the three paths outlined in [B-M] mentioned above, but including now all the intersections of quadrics and not only the moment-angle manifolds:

\bigskip

1) In Section 1 we identify very general families of manifolds that are indeed
diffeomorphic to connected sums of sphere products. The main result there is Theorem 1.3 (of which the following is a simplified version), where $Z$ denotes
an intersection of quadrics and $Z^J$ is an infinite family obtained from $Z$ by increasing the dimension of its coordinate spaces:

\bigskip

\textbf{Theorem.} \textit{Assume $Z$ is of dimension $2c$ and
$(c-1)$-connected where $c\ge 2$. Then any $Z^J$ of dimension at
least 5 is diffeomorphic to a connected sum of sphere products.}

\bigskip

Particular members of such families are those moment-angle manifolds for which the result answers the Bosio-Meersseman conjecture in [B-M]. The dimension hypothesis only excludes the simply connected 4-dimensional ones, for which the result should also be true but our methods only show that they are  \textit{homeomorphic} to a connected sum of copies of $S^2\times S^2$. A weaker result in the odd dimensional case is given in Theorem 1.4, while Theorems 1.1 and 1.2 are relative results used in the proofs of the main theorems, but which have other applications. As a byproduct of the proofs, a simpler and neater proof of the result in [LdM1] is also obtained. 

\bigskip

2) In section 2 we give an explicit topological description of some of the transitions, mainly that of cutting off a vertex or an edge of the associated polytope (operations $Z_\textsf{v}$ and $Z_\textsf{e}$). We show that, under simple assumptions, these operations preserve connected sum of sphere products and adequately modified (operations $Z_\textsf{v}\, '$ and $Z_\textsf{e}\, '$) they can be combined with the theorems in section 1 to give new infinite families that are diffeomorphic to connected sums of sphere products, the final result along this line being:

\bigskip

\textbf{Theorem 2.4:} \textit{If $Z = Z(P)$ is a connected sum of sphere products and is simply connected of dimension at least 5, then any manifold obtained from $Z$  by repeatedly applying the $Z^J$, $Z_\textsf{v}\, '$ and $Z_\textsf{e}\, '$ constructions (any number of times each and in any order) is also a connected sum of sphere products.}

\bigskip

The initial manifold can be any sphere, any highly connected manifold $Z(P)$ as in Theorems 1.3 and 1.4 or any simply connected manifold with $m - d = 3$ of dimension at least 5 that is not a triple sphere product ([LdM1]).

\bigskip

The first results in section 2 include cases where the initial manifold is quite general and we make a similar conjecture specific for moment-angle manifolds. A proof of that conjecture under a restrictive hypothesis is enough for us to describe the topology of other important examples taken from [B-M]: the manifolds associated to the truncated cube.

\bigskip

3) In section 3 we use our knowledge of the above examples to give a proof that their cohomology rings are not isomorphic as ungraded rings, not even when taken with $\Z_2$ coefficients. Thus we show that the product rule for the moment-angle manifolds has to be drastically modified in the general situation. This contradicts a result by M. de Longueville ([dL]). We only state the modified rule, and leave the details for another publication ([G-LdM]).

\bigskip

In section 0 we recall the necessary definitions and known results and in section 2.1 some elementary topological constructions are defined and explored. In the appendix we state and prove some results about specific differentiable manifolds, which are used in sections 1 and 2.

\bigskip

In the process of solving these problems several new questions and conjectures have arisen. Some extensions of our results are immediate, others should be possible by continued work along the same lines. Others seem to need a whole new approach.

\subsection*{Acknowledgements.}

This research was initiated by a discussion with Fr\'ed\'eric
Bosio and Laurent Meersseman around the article [A-LL] and continued
through fruitful conversations with them and with Tony Bahri, Martin
Bendersky, Fred Cohen, Vinicio G\'omez Guti\'errez, Francisco Gonz\'alez Acu\~na
and Alberto Verjovsky.

\bigskip

\noindent Part of this research was done while the second author enjoyed the
warm and stimulating hospitality of the NYU Courant Institute, for
which he is most grateful to Sylvain Cappell. Invitations to
Princeton University by William Browder and to the CUNY Graduate
Center by Martin Bendersky and Maxim Laurentiu were also very useful
and encouraging. He was also partially supported by DGAPA-UNAM grant PAPIIT-IN102009.

\newpage

\section* {0. Preliminaries.}

\subsection {Basic definitions.}

Let $m>d>0$ be integers and $k=m-d-1$,
$$\Lambda_i\in \R^k, i=1,\dots,m, \,\, \Lambda=(\Lambda_1, \dots, \Lambda_m)$$

If $I\subset \{1,\dots,m\}$ we will denote by $\Lambda_I$ the
sequence of $\Lambda_i$ with $i\in I$.

\vspace{.15in}

We will assume the (generic) \textit{weak hyperbolicity} condition:
\vspace{.25in}

\emph{If $I$ has $k$ or less elements then the origin is not in the
convex hull of $\Lambda_I.$}

\bigskip

Let $Z=Z(\Lambda)$ be the submanifold\footnote{ These intersections of quadrics were first
studied in [W3], [Ch], [H] and [LdM1]. The manifolds $Z^{\C}$ were first studied in [C-K-P], [Ch] and [LdM2] where they arose (and reappeared recently, see [Ch-LdM]) in the context of dynamical systems.

$Z$, $Z^{\C}$ are essentially the same as the manifolds constructed independently in [D-J], known as
\textit{universal abelian covers} of polytopes and  \emph{moment-angle manifolds}, respectively. $Z$ and the constructions $Z_+$, $Z'$, $Z^J$ described below, can be
expressed in the framework of the \textit{generalized moment-angle complexes} defined in [B-B-C-G].}
of $\R^m$ given by the equations

$$
\mathop{\Sigma}\limits_{i=1}^n \Lambda_i x_i^2=0
$$

$$
\mathop{\Sigma}\limits_{i=1}^n x_i^2=1.
$$

It is immediate that this manifold is non-empty if, and only if, the origin of $\R^k$ is in the convex hull of the $\Lambda_i$.

\vspace{.15in}

$Z$ is smooth, since the weak hyperbolicity condition is equivalent to the regularity of the system of its defining equations.
It is immediate that it is stably parallelizable. A simple construction (see section 1) shows that it is always the boundary of a parallelizable manifold.

\vspace{.15in}

Let $Z^{\C}=Z^{\C}(\Lambda)$ be the manifold defined in $\C^n$ by the equations:

$$
\mathop{\Sigma}\limits_{i=1}^n \Lambda_i \left|z_i\right|^2=0
$$

$$
\mathop{\Sigma}\limits_{i=1}^n \left|z_i\right|^2=1.
$$

Of course, $Z^{\C}$ is a particular (but very important) case of a $Z$ manifold if we write its
equations in real coordinates in $\R^{2n}$, each $\Lambda_i\in R^k$ appearing twice.

\vspace{.15in}

It was shown in [LdM1] that for $k=2$, all manifolds $Z^{\C}$ and almost all simply connected manifolds $Z$
are either empty, a product of three spheres or a connected sum of manifolds, each of
which is a product of two spheres. We will call this type of sums
\textit{connected sum of sphere products} (implicitly implying that
each summand is the product of only two spheres).

It is easy to see that this collection of manifolds is closed under products so one gets examples for high $k$ where $Z$ is a product of any number of spheres and any number of connected sums of sphere products, which have therefore free homology groups and easily described cohomology rings.

In [B-M] it was shown that for $k\ge 3$ the manifold $Z^{\C}$ may be considerably more complicated (see also [Ba2]) and that for $k\ge 5$ its homology may have a lot of torsion. But they also conjectured that under certain conditions one obtains again the same kind of connected sums. The results of this article imply that this Bosio-Meersseman conjecture is true and we obtain the same conclusion in many more cases.

\subsection {Group actions and polytopes.}

The manifold $Z$ admits a $\Z_2^m$ action by changing the signs
of the coordinates. The quotient is a simple polytope $P$ which can
be identified with the intersection of $Z$ and the first orthant of
$\R^m$. It follows that $Z$ can be reconstructed from this intersection by
reflecting it on all the coordinate hyperplanes.

By a simple change of coordinates $r_i=x_i^2$, this quotient can be identified with the $d$-dimensional convex polytope $P$ given by
$$
\mathop{\Sigma}\limits_{i=1}^n \Lambda_i r_i=0.
$$

$$
\mathop{\Sigma}\limits_{i=1}^n r_i=1.
$$

$$r_i\ge 0.$$

The weak hyperbolicity condition is equivalent to the fact that $P$ is a \textit{simple polytope}, meaning that each vertex is exactly in $d$ facets of $P$.

The \textit{facets} of $P$ (i.e., its $(d-1)$-dimensional faces) are the non-empty intersections of $P$ with the coordinate
hyperplanes $\{r_i=0\}$. It follows that a non-empty $Z$ is connected if, and only if, all these intersections are non-empty, or equivalently, if $P$ has exactly $m$ facets, since otherwise different components of $P$ would lie in different sides of an hyperplane $\{x_i=0\}$. In that case $Z=Z_* \times \Z_2^h$ where $Z_*$ is a connected
component of $Z$ and $h$ is the number of those coordinate hyperplanes not touching $Z$.

If we introduce a new facet by intersecting $P$ with a half space $H=\{\mathop{\Sigma}\limits_{i=1}^n a_i r_i \ge a\}$ whose boundary does not contain any vertex of $P$ we obtain a new simple polytope $P_H$. By introducing a new coordinate $r_0= \Sigma_{i=1}^n a_i r_i- a$ we can put $P_H$ in $\R^{n+1}$ with $r_0\ge 0$. By simple manipulations of the equations
of $P_H$ and changes of coordinates of the type $r_i'= b_i r_i$, $i=0, \dots, n$ we get a system of equations for $P_H$ of the same form as that of $P$.

This implies, by induction on the number of facets, that any simple polytope can be expressed in the above form and is therefore the quotient polytope of a connected manifold $Z$.\footnote{ This \emph{Realization Theorem} is alluded to in [W3] as a general correspondence known as the \emph{Gale diagram} of the polytope. The \emph{perfect} correspondence between polytopes and the manifolds $Z$ is mentioned in [LdM], but its precise statement and proof were not included after the author learned about the Gale diagram from Wall's paper. A detailed proof using the Gale Transform is given in [B-M], Theorem 0.14. The simple proof by induction on the number of facets of $P$ is due to Vinicio G\'omez Guti\'errez.}

This means that the connected $Z$ coincide with the manifolds constructed in [D-J], called by them \textit{universal abelian covers}.

We can thus refer to the manifold $Z(P)$, by which we mean the diffeomorphism type of any connected manifold $Z(\Lambda)$
whose quotient polytope is $P$.\footnote{One can deduce from [B-M] that combinatorially equivalent polytopes are in fact diffeomorphic as
manifolds with corners, and thus there is really only one diffeomorphism type $Z(P)$. We will not be using this rigidity result since all our constructions
can be given in terms of the equations and the uniqueness follows in general from the results themselves.}

It is known that $Z$ is $(c-1)$-connected if, and only if, the intersection of any $c$
facets of $P$ is non-empty.\footnote{ The homology groups of $Z$ follow from the computations given in [LdM1] or [B-B-C-G].
An argument for the fundamental group is given in [LdM1].} Following [B-M],  such a $P$ is called \emph{dual $c$-neighborly} and simply \emph{dual neighborly} if it is $[d/2]$-dual neighborly. The B-M conjecture says that if $P$ is even-dimensional and dual-neighborly then $Z^{\C}$ is a connected sum of sphere products.

\bigskip

The manifold $Z^{\C}$ admits an action of the $m$-torus $T^m=(S^1)^m$ by multiplication of unit complex numbers
on each of the coordinates. The quotient by this action is the same simple polytope $P$ and $Z^{\C}$ can be
reconstructed from $P\times T^m$ by making the obvious identifications on the facets of $P$. The orbits of this
action over a facet $F$ correspond to the points where some coordinate $z_i$ vanishes and are thus the points where
the corresponding factor of the torus acts trivially.

If $P$ has exactly $m$ facets, any orbit of the action of $T^m$ on $Z^{\C}$ is null-homotopic in $Z^{\C}$.
This is because it can always be deformed to an orbit over a facet of $P$, where one of the factors of $T^m$
acts trivially and so the orbit collapses into an orbit of a $T^{m-1}$-action. This new orbit can be deformed to a $T^{m-1}$
orbit over another facet of $P$ where another factor acts trivially. Since any factor of $T^m$ acts trivially on some facet we
end up with a deformation of the original orbit into a point.

This proves that in that case $Z^{\C}$ is simply connected, since any loop in it can be moved away from the manifolds lying over the facets of $P$ (which have codimension 2) and so lies in the open set $interior(P)\times T^m$ and can therefore be deformed
into a single orbit which in turn deforms into a point. From the same homology computations as above, $Z^{\C}$ is then
$2c$-connected if, and only if, $P$ is dual $c$-neighborly, and therefore always actually 2-connected.

On the contrary, if $P$ has $m-h$ facets then  $Z^{\C}=Z^{\C}_* \times \T^h$. (For a proof see [B-M], Lemma 0.10). In that case the orbits of $T^m$ are not all null-homotopic and $Z^{\C}$ is not simply connected.

The 2-connected $Z^{\C}$ coincide with the manifolds constructed in [D-J], called now \textit{moment-angle manifolds}.
This fact follows from the realization of a polytope as the quotient of a connected manifold $Z$ (see footnote 3).

We can therefore refer to the manifold $Z^{\C}(P)$, by which we mean the diffeomorphism type of any 2-connected manifold $Z^{\C}(\Lambda)$
whose quotient polytope is $P$. It was shown in [B-M] that this diffeomorphism type and the $T^m$ action are determined by the combinatorial type of $P$.

\newpage

\section {Opening books.}

Let $\Lambda'$ be obtained from $\Lambda$ by adding an extra
$\Lambda_1$ (so $\Lambda'$ satisfies also weak hyperbolicity) which we interpret as the coefficient of a new 
variable $x_0$, so we get the manifold $Z'$:

$$
\Lambda_1 (x_0^2 + x_1^2)+\mathop{\Sigma}\limits_{i>1} \Lambda_i x_i^2=0
$$

$$
x_0^2 + x_1^2 + \mathop{\Sigma}\limits_{i>1} x_i^2=1.
$$

Let $Z_+$ be the intersection of $Z$ with the half space $x_1\ge 0$
and $Z'_+$ the intersection of $Z'$ with the half space $x_0\ge 0$.
The boundary of $Z'_+$ is $Z$. This shows that $Z$ is always the
boundary of a parallelizable manifold.

$Z_+$ admits an action of $\Z_2^{n-1}$ by changing signs on all the variables except $x_1$ and the quotient is again $P$.
In other words, $Z_+$ can be obtained from $P\times \Z_2^{n-1}$ by making all the identifications on the facets of $P$ except
for the one corresponding to $x_1=0$ (or, equivalently, by reflecting $P$ on all the coordinate hyperplanes except $x_1=0$).

Consider also the manifold $Z_0$ which is the intersection of $Z$
with the subspace $x_1 = 0$. $Z_0$ is the boundary of $Z_+$. So
$$Z_0\subset Z_+ \subset Z \subset Z'_+ \subset Z'$$

For example,  $Z_+(\Delta^n)=D^n$ that can be obtained from $\Delta^n\times \Z_2^n$ by making the identifications on all the faces
of $\Delta^n$ but one (or, equivalently, by reflecting the spherical simplex $S^n \bigcap \R_+^{n+1}$ on all the coordinate
hyperplanes but one). Also $Z_+(\Delta^n \times I)$ can be either $D^n\times S^1$ if we do not reflect it on a facet $\Delta^{n-1} \times I$
or $S^n\times D^1$ if we do not reflect it on a facet $\Delta^n \times {pt}$.

The quotient polytope $P'$ of $Z'$ is the \emph{book} on $P$ obtained from $P\times I$ by collapsing $P_0\times I$ into $P_0$
where the latter is the facet of $P$ corresponding to $x_1=0$ (that is, $P_0$ is the quotient polytope of $Z_0$).\footnote{ This construction is implicit in [LdM1] and was much used in the process of proving its main Theorem 2. The picture of
the pentagonal book appears already in [W3], p.413. The name \emph{book} was used in [Me] and [B-M], while in the literature stemming
from [D-J] it is known as the Buchstaber construction, see [B-P].}

\bigskip

$S^1$ acts on $Z'$ (rotating the coordinates $(x_0,x_1)$) with fixed
set $Z_0$. Its quotient can be identified with $Z_+$. The map
$$(x_0,x_1,x_2,\dots,x_m)\mapsto (\sqrt{x_0^2+x_1^2},x_2,\dots, x_m)$$
is a retraction from $Z'$ to $Z_+$ which restricts to the retraction
from $Z$ to $Z_+$
$$(x_1,x_2,\dots , x_m)\mapsto (\left|x_1\right|,x_2,\dots , x_m)$$

\begin{center}
\includegraphics[width=2.8in]{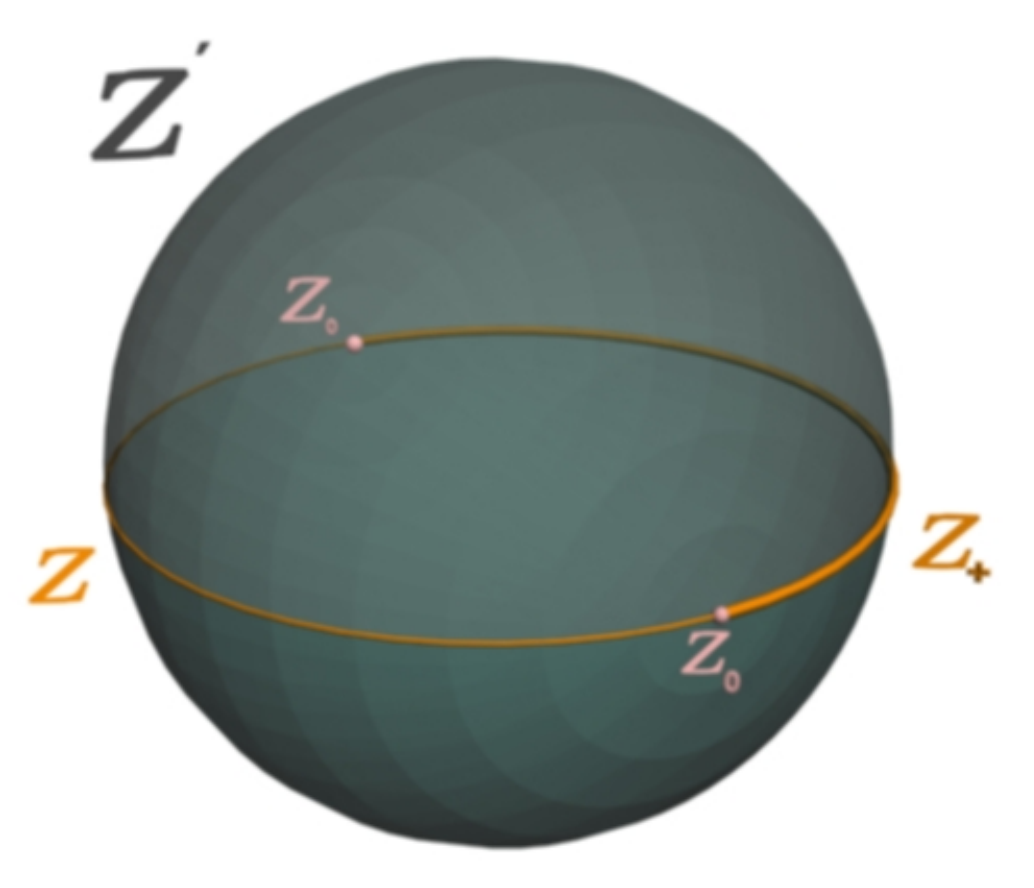}
\end{center}

\bigskip

\bigskip

\bigskip

Observe further that this retraction restricted to $Z'_+$ is
homotopic to the identity: the homotopy preserves the coordinates
$x_i,i\ge2$ and folds gradually the half space $x_0\ge 0$ of the
$x_0,x_1$ plane into the ray $x_0=0,x_1\ge 0$ preserving the
distance to the origin:

\bigskip

\bigskip

\bigskip

\begin{center}
\includegraphics[width=4in]{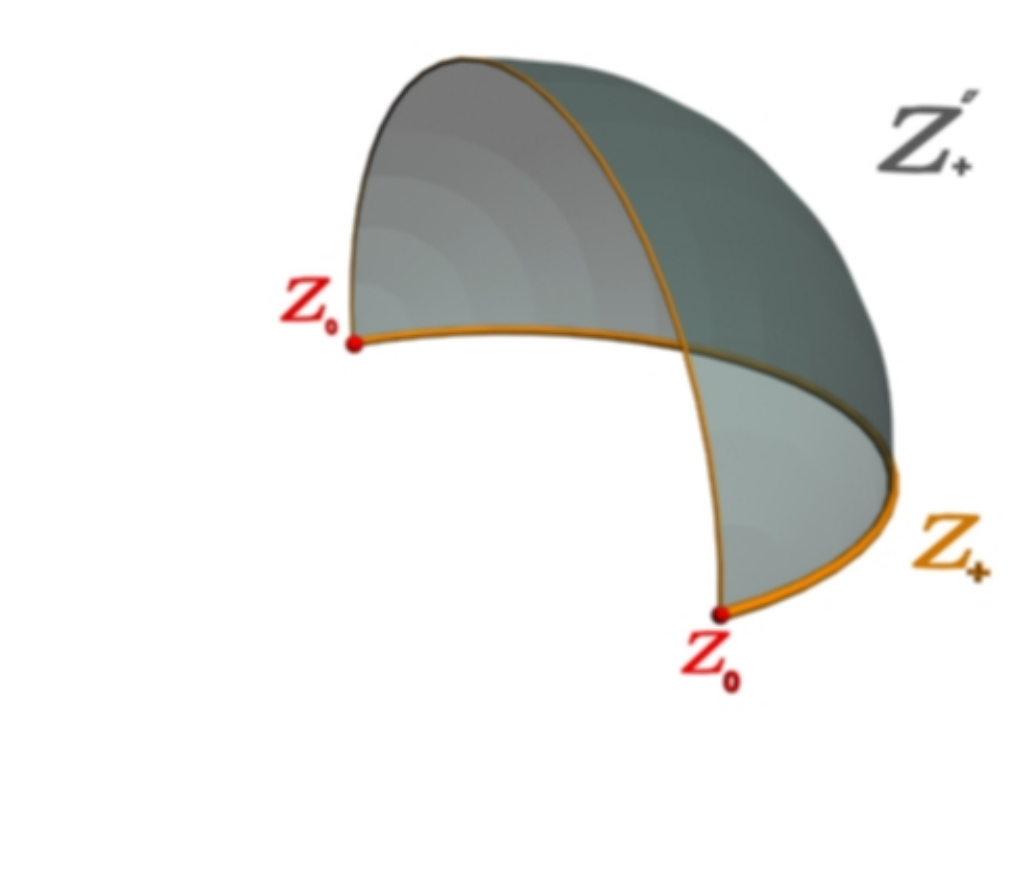}
\end{center}

\bigskip

So $Z$ is the double of $Z_+$ and $Z'$ is the double of $Z'_+$, and $Z'$ is the open book with page $Z_+$ and trivial holonomy.\footnote{In the notation of section 2.1,
$ Z= \mathcal{D}(Z_+), Z'= \mathcal{D}(Z'_+)$ and $ Z'= \mathcal{TOB}(Z_+)$.}

The inclusion of $Z=\partial Z'_+$ into $Z'_+$ is homotopic to the retraction from $Z$ to $Z_+$.

The construction of $Z'$ from $Z$ can be generalized as follows:
Let $J=(j_1,\dots,j_m)$ be a vector of positive integers and $\vert J\vert=\sum j_i$. Then
we construct a new configuration $\Lambda^J$ obtained from $\Lambda$
by repeating $j_i$ times the coefficient $\Lambda_i$ for each $i$ and
thus a new corresponding manifold $Z^J$. The manifold $Z^J$ has a natural $O(j_1) \times \dots \times O(j_n)$ action whose quotient is the same simple polytope $P$, as well as its own $\Z_2^{\vert J\vert}$ action with quotient a polytope $P^J$.

For example, if $J=(2,1,\dots,1)$ we have $Z^J=Z'$. If $j$ is a
positive integer we denote by $\underline{j}$ the vector
$(j,j,\dots,j)$. Thus $Z^{\underline{1}}=Z$ and $Z^{\underline{2}}=Z^\C$.

Independently of any construction related to quadrics, the $Z^J$ construction appeared in [B-B-C-G2] (already in its first draft of 2008) in the context of toric manifolds and generalized moment-angle complexes. There it is shown how to describe explicitly the polytope $P^J$ and how to construct a new toric manifold $M^J$ associated to $P^J$ from a toric manifold $M$ associated to $P$.

In fact, all our constructions can be expressed within the context of the generalized moment-angle complexes of [B-B-C-G] (but not inside any of the many previous generalizations of moment-angle manifolds) and the general splitting formula proved there provides an excellent geometric tool for understanding (among other things) the relations among the homology groups of all these manifolds. It is out of this understanding that were discovered the geometric relations used in the following proofs.

We will prove that certain families of manifolds $Z$ are connected sums of sphere products.
The rough idea of the proof in each case is as follows:

\bigskip

(i) We show that the lower dimensional manifolds in the family, being highly connected,
are connected sums of sphere products.

\bigskip

(ii) Then we use the geometric relation between the manifolds $Z$ and $Z'$ to prove the
induction step that allows us to cover all the high dimensional ones.

\bigskip

Part (i) can be derived from the classification results about highly connected manifolds
(see, among others, [Sm], [K-M], [W1], [W2]). This type of result usually involves
practically all the techniques of simply-connected surgery ([Br], [K-M]) which include the number-theoretic and algebraic analysis of the middle dimensional quadratic forms in order to represent homology classes by embedded spheres with trivial normal bundle, the construction and simplification of an adequate cobordism and the use of the $h$-cobordism theorem.

In our case we will not use those results but will use instead a variation of the arguments used in the induction step (ii) to prove (i) directly using the simple nature of our manifolds for which those constructions are elementary.

\bigskip

Part (ii) follows the spirit of those classification theorems, but is considerably simpler. We use the obvious fact that if $X$ is a simply-connected connected sum of sphere products, then a basis of each homology group below the top dimension can
be represented by embedded spheres with trivial normal bundle. A simple geometric lemma
lets us represent various combinations of them by embedded spheres with trivial normal bundle inside $X\times I$ and disjoint from each other.

Then we use the natural cobordism $Z'_+$ which has the crucial but simple property that
all its homology comes from its boundary $Z$. Combining this with the previous fact we
can use Theorem A1 in the appendix to show that $Z'_+$ is a connected sum along the boundary of trivial bundles over spheres, so the induction step follows easily for its double $Z'$.

We start by the induction step (Theorem 1.1) and its variation that allows us to simplify the starting point of the induction (Theorem 1.2). These are interesting by themselves since they can be used in other situations. Then we will prove our two main results (Theorems 1.3 and 1.4).

\bigskip

\textbf{Theorem 1.1.} \textit{Assume $Z$ is simply connected and of
dimension $d\ge 5$. If $Z$ is a connected sum of sphere products
then $Z'$ is also a connected sum of sphere products.}

\bigskip

\textbf{Proof:}

\bigskip

We start with a

\bigskip

\textbf{Lemma.} \textit{If $X$ is a manifold of dimension $d$ and
$\alpha_1,\alpha_2 \in H_i(X)$ can be represented by disjoint
embedded spheres $S_1,S_2$ with trivial normal bundle, then if $i\le
d-2$, $\alpha_1+ \alpha_2$ can be represented by an embedded sphere
with trivial normal bundle and disjoint from $S_1$ and $S_2$.}

\bigskip

\textit{Proof of the lemma}: Taking disjoint tubular neighborhoods
of $S_1,S_2$ and joining them with a thin tube, we obtain a new
manifold $X'$ diffeomorphic to the connected sum along the boundary
of two copies of $S^i \times D^{d-i}$. It is enough to prove the
lemma for $X'$. We can take for $X'$ a standard model constructed as
follows: take $S'_1, S'_2$ to be round spheres in $\R^{i+1}\subset
\R^d$ of radius $1/2$ centered in $(1,0,\dots,0)$ and
$(-1,0,\dots,0)$ together with the straight line segment $\sigma$
joining $(-1/2,0,\dots,0)$ and $(1/2,0,\dots,0)$. Now consider an
$\epsilon$ neighborhood of this configuration in $\R^d$.

One can assume that the spheres with their standard orientations in
$\R^{i+1}$ correspond to the classes $\alpha_1,\alpha_2 \in H_i(X)$:
Because if $S'_i$ corresponded to $-\alpha_i$, one can reflect
the product $S^i \times D^{d-i}$ on the hyperplane $x_2=0$ of $\R^d$
and twist the tube around $\sigma$ accordingly.\footnote{ When
$i=n-1$ this last step cannot be achieved. For example, in
$X=S^{d-1}\times S^1$ $\alpha + \alpha$ (for $\alpha$ the generator
of $H_{d-1}(X)$) is not representable by an embedded sphere even if
we can find two disjoint spheres each representing $\alpha$.}

Taking the intersection of this neighborhood with $\R^{i+1}$ we
obtain a manifold whose boundary is formed by three embedded spheres
with trivial normal bundle in $\R^{i+1}$ and therefore also in
$\R^d$. The outside sphere (adequately oriented) represents the sum $\alpha_1+ \alpha_2$
and is disjoint from $S_1$ and $S_2$. \textit{QED}

To apply the lemma to the situation of the theorem, observe that,
since $Z$ is a connected sum of sphere products, all its homology
groups are free and all the generators of the homology below the top
dimension are represented by embedded spheres with trivial normal
bundle. Furthermore, they can all be made disjoint inside $Z\times
I$ by displacing them to different levels of the coordinate $t\in
I$. They all have codimension at least $2$ in $Z$ since $Z$ is
simply connected, so they have codimension at least 3 in $Z\times
I$. The lemma implies that not only the generators, but all elements
of the homology below the top dimension can be represented by
embedded spheres in $Z\times I$ with trivial normal bundle\footnote
{This may not be true inside $Z$ itself: in $S^{2c}\times S^{2c}$
the sum of the generators can be represented by the diagonal sphere,
but its normal bundle is the tangent bundle of $S^{2c}$ which not
trivial. No other embedded sphere representing that sum can have
trivial normal bundle since its self-intersection number is 2.}
where any finite collection of them can be represented by mutually
disjoint ones.

As we showed early in this section, the inclusion of $Z=\partial Z'_+$
into $Z'_+$ is homotopic to a retraction. It therefore induces a
split epimorphism in homology\footnote{This was proved for $k=2$ in
[LdM1] using an explicit computation of the homology groups involved
in that special case. A first proof in the general case used the stable
splitting in [B-B-C-G] and the relationship between the polytopes $P$ and $P'$, but its generality suggested the search for a simpler geometrical argument.}
and in the fundamental group.

So $Z'_+$ is simply connected and its homology is also free and is trivial
in dimensions $d-1$ and $d$.

A basis for the homology of $Z'_+$ can then be represented
by mutually disjoint embedded spheres with trivial normal bundle
lying in a collar neighborhood of its boundary $Z$. They all have
codimension at least 3, so the complement of their union in $Z'_+$
is still simply connected.

So all the hypotheses of Theorem A1 in the appendix are fulfilled and we conclude
that $Z'_+$ is the connected sum along the boundary of products of the
form $S^p \times D^{d+1-p}$

Finally, $Z'$, being the double of $Z'_+$, is the connected sum of
manifolds of the form $S^p \times S^{d+1-p}$. \textbf{QED}.

\bigskip

Observe that all we have used about $Z$ is that it is simply
connected, all its homology groups are free and every homology class
below the top dimension can be represented by an embedded sphere
\textit{inside $Z\times I$} with trivial normal bundle. And the
proof gives that $Z$ is in fact a connected sum of sphere products.
So, including the full induction on the size of $J$, we have
actually proved more:

\bigskip

\textbf{Theorem 1.2.} \textit{Assume that $Z$ is simply connected and
of dimension $d\ge 5$ and that all its homology groups are free and
every homology class below the top dimension can be represented by
an embedded sphere inside $Z\times I$ with trivial normal bundle.
Then $Z$ is a connected sum of sphere products and so is $Z^J$ for
any $J$.}

\bigskip

As a first application, assume $Z$ is of dimension $2c$  and
$(c-1)$-connected (which means that the corresponding polytope is
\textit{dual neighborly}, see section 0.2).
If $c\ge 3$ then $Z$ clearly satisfies the hypotheses of Theorem 1.2
and every $Z^J$ is a connected sum of sphere products. For $c=2$ we
only get that $Z$ is $h$-cobordant to a sum of copies of $S^2\times S^2$.
But we can apply the construction in the proof of the theorems one more
step up: one obtains that $Z'$ is the boundary of a manifold $Z_+''$
which deforms down to $Z'_+$. In $Z'_+$ all the homology classes
are in dimension 2 and can be represented by embedded spheres with
trivial normal bundle, so the same is true for $Z_+''$. Now the dimension of $Z_+''$ is 6 and now Theorem A1 shows again that it is a connected sum along the boundary of copies of $S^2\times D^4$ and so its boundary $Z'$ is a connected sum along
the boundary of copies of $S^2\times S^3$. Applying now Theorem 1.1 inductively we obtain that all higher $Z^J$ are connected sums of products of spheres. So we have:

\bigskip

\textbf{Theorem 1.3.} \textit{Assume $Z$ is of dimension $2c$ and
$(c-1)$-connected. Then:}

\textit{1) If $c \ge 3$, $Z$ and any $Z^J$ are connected sums of sphere products.}

\textit{2) If $c=2$, any $Z^J$ of dimension at least 5 is a connected sum of sphere products.}

\bigskip
\textbf{Remarks.} 1. This includes the Bosio-Meersseman conjecture ([B-M])\footnote{ Thus
the counterexample to this conjecture announced in [A-LL] is no
such. To see directly that the argument there is not correct one can
use the epimorphism part of the Hurewicz theorem to check that
$\pi_6(Z)$ is not as stated.} for $c \ge 2$ by taking the special case $J=\underline{2}$.

\bigskip

2.- The result 1) should be true for any $c\ge 1$ and without any restrictions. For $c=2$ the only missing case is $Z$ itself, where
the topological h-cobordism theorem in dimension 4 shows that it is \textit{homeomorphic} to a connected sum of copies of
$S^2\times S^2$. For $c=1$, $Z$ has dimension 2 and is trivially a connected sum of sphere products, and so is $Z'$ by Theorem 2.1
since $P'$ can be obtained from a simplex by cutting vertices. For the case of polygons with at most 8 sides it can be proved using the techniques of section 2. For $J$ at least $\underline{2}$ it would follow from our conjecture preceding Theorem 2.2.
This would also follow from Theorem 6.3 of [B-M] but the proof of this result is not given and it is not clear to us how such a proof would avoid the problems
we find in proving the above mentioned conjecture. Also, the result from [Mc] alluded to in [B-M] is proved only in the PL-category.

\bigskip

3.- Our result provides also a new, neater proof of the result in [LdM1], starting from the basic $2k$-dimensional $(k-1)$-connected manifolds when $k>2$ and with the simply-connected 5-dimensional ones when $k=1,2$, all of which can be covered by the theorems of this section.

\bigskip

In the odd dimensional case we have for the moment only weaker results:
If $Z$ is of dimension $2c+1$  and $(c-1)$-connected for $c>1$ we need
to assume that all its homology groups are free, a condition that is automatic in the even dimensional case. Nevertheless, it should be mentioned that the example of
a manifold $Z$ with torsion in [B-M] does not have the high connectivity assumed here, so that hypothesis might still be redundant.

On the other hand we have a slightly weaker conclusion because we
cannot exclude the possibility of having a connected sum of sphere
bundles over spheres which are stably trivial over the embedded
$(c+1)$-spheres but not trivial. When $c$ is odd the total space of those
bundles has torsion so this possibility is excluded by our hypothesis.
When $c$ is even the only possibility is the unit tangent bundle of
of $S^{c+1}$ which is trivial only for $c=2,6$. In this case they
cannot be discarded by their torsion. Yet it is hard to see how those
non-trivial bundles could have the high $\Z_2^m$ symmetry of our
manifolds.

In any case this problem disappears for the manifold $Z'$ by going
up to the manifold $Z''_+$ as above\footnote{ Otherwise, remark that
the double of a connected sum along the boundary of stably trivial
$D^{c+1}$-bundles over $S^{c+1}$ is a connected sum of trivial
$S^{c+1}$-bundles over $S^{c+1}$.}:

\bigskip

\textbf{Theorem 1.4.} \textit{Assume $Z$ is of dimension $2c+1$ and
$(c-1)$-connected where $c\ge 2$ and assume further that all the
homology groups of $Z$ are free. Then:}

(\textit{i) $Z$ is a connected sum of stably trivial sphere bundles over
spheres and if $J \ne\underline{1}$ then $Z^J$ is a connected sum of
sphere products.}

\textit{(ii)  If $c$ is odd or if $c=2,6$, then $Z$ is a connected sum
of sphere products.}

\bigskip

The case $J=\underline{2}$ is a weak version of the Bosio-Meersseman
conjecture in the odd-dimensional case.

Theorem 1.4 is definitely not true for $c=1$: for the cube, $Z=S^1 \times S^1 \times S^1$.

\bigskip

In Theorems 1.3 and 1.4 the actual dimensions of the spheres can be deduced
from the computation of the homology groups of the manifolds. These can be
deduced from the combinatorics of the polytope as in [LdM1], [B-P] or [B-M],
or by using the stable splitting proved in [B-B-C-G]. Observe that we
have also proved that the halved manifolds $Z'_+$ are connected sums
along the boundary.

\bigskip

A general proof that does not use the $h$-cobordism theorem and that
allows one to actually \textit{see} the spheres is still badly
needed.

\newpage

\section {Cutting corners and making ends meet.}

We study now the effect on the manifold $Z$ of cutting off parts of the polytope
$P$ of dimension $d$. The general idea can be seen in the following situation:

\begin{center}
\includegraphics[width=9in]{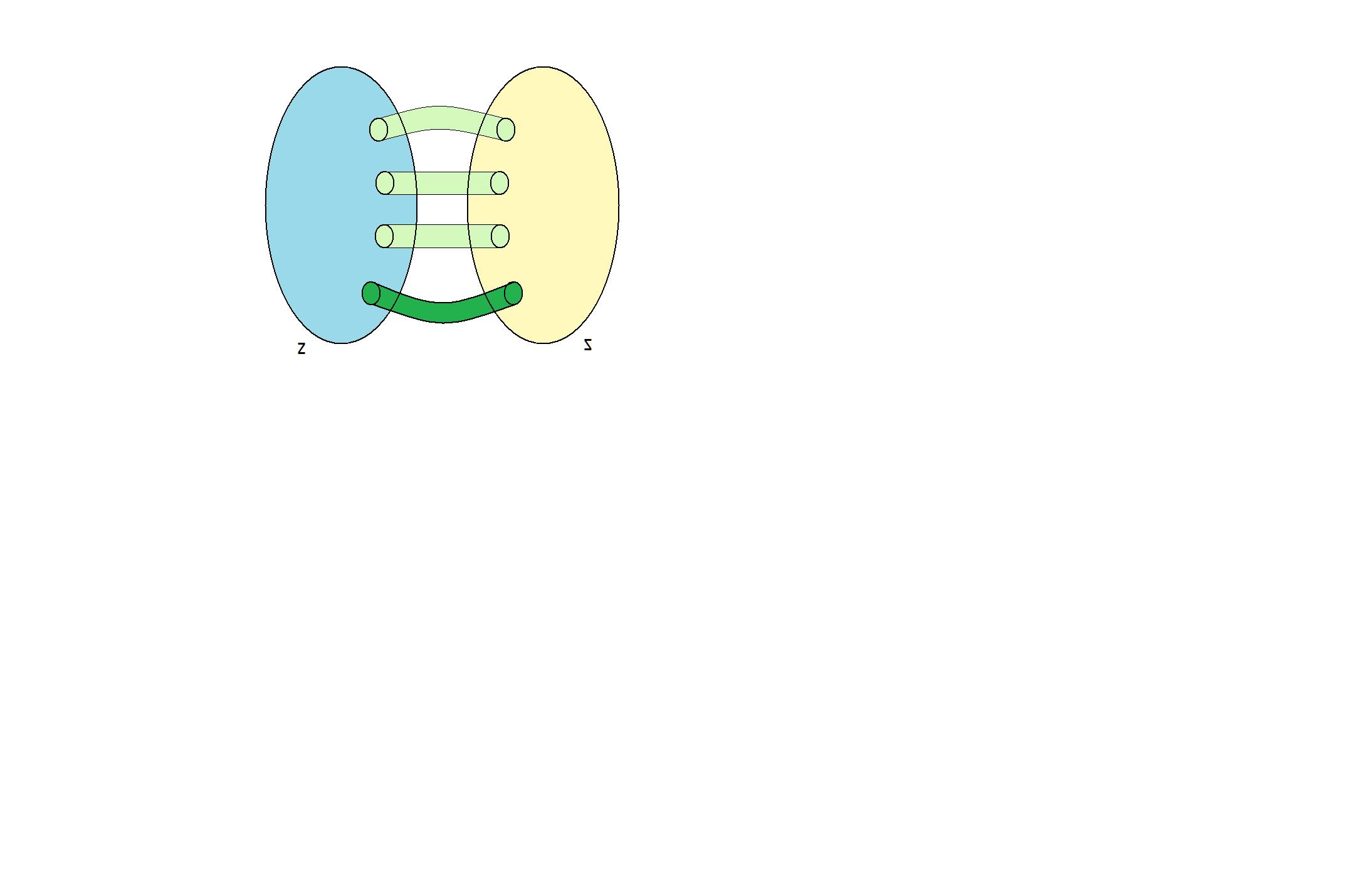}
\end{center}

\vskip -19.55em

Cutting off a vertex of $P$ gives a new polytope $P_\textsf{v}$ and its corresponding manifold $Z\textsf{v}$. In $Z$ this has the effect of removing a
number of open balls thus obtaining a manifold $Z_{\textsf{v}+}$ with boundary. To obtain the manifold $Z\textsf{v}$ we have to take the \textit{double} of $Z_{\textsf{v}+}$ by taking two copies of it and identifying their boundaries (consisting each of a number of spheres). In the figure we have simulated this by joining those spheres by tubes $D^1 \times S^{d-1}$.

We obtain the same result by taking two copies of $Z$ and performing several times the operation of removing a disk from each copy and joining their boundary spheres by a tube. Making the operation the first time amounts to taking the connected sum of the manifolds (imagine only the lower two disks and the darker tube). Since the resulting manifold is connected, each new time we perform the operation amounts to attaching a handle or, equivalently, taking the connected sum with $S^1\times S^{d-1}$.
Thus the final result is the connected sum of two copies of $Z$ and a number of copies of $S^1\times S^{d-1}$ (one less than the number of open balls removed from $Z$).

To make precise this construction and some of its variants we introduce now some simple topological constructions and explore the relations among them.

\subsection {Doubles, Open Books and Connected sums.}

We work in the smooth category. In what follows we consider many manifolds with corners and we will use implicitly and systematically the fact that they have natural smoothings.

\bigskip

If $Q$ is a manifold with boundary, we denote its boundary by $\partial Q$ as usual.

\bigskip

Let $\mathcal{D} Q$ be the \textit{double} of $Q$:

\bigskip

$$\mathcal{D} Q =\partial(Q \times D^1).$$

Let $\mathcal{TOB}(Q)$ be the open book with trivial holonomy with leaf $Q$ and binding $\partial Q$, that is

$$\mathcal{TOB}(Q)=Q\times S^1 \cup \partial Q\times D^2=\partial(Q\times D^2).$$

\bigskip

If $M$ is a manifold without boundary and $M_{-1}$ is $M$
minus an open ball, $\mathcal{D} M_{-1}$ is diffeomorphic to
$M\#(-M)$. In all cases we will consider below, $M$ admits an orientation
reversing diffeomorphism, so $\mathcal{D} M_{-1}$ will be diffeomorphic to $M\# M$.

\bigskip

Let $\mathcal{G}M=\partial(M_{-1} \times D^2)$.\footnote{ $\mathcal{G}$ stands for \textit{gyration} and also
for Gonz\'alez Acu\~na (Fico) who defined and studied this operation in [GA].}

\bigskip

So $$\mathcal{G}M=\mathcal{TOB}(M_{-1}).$$

\bigskip

As a simple example, $$\mathcal{G}(S^n)=S^{n+1}.$$

\bigskip

If $Q,Q'$ are connected $d$-manifolds, then we will consider their connected sum $Q\#Q'$ (connecting the \emph{interiors} of them).
If $Q,Q'$ have also connected non-empty boundaries we can take their connected sum along the boundary $Q\amalg Q'$.

\bigskip

We explore now the relations between all these operations:

\bigskip

Clearly, if $M$ and $N$ are closed and connected and $Q,Q'$ are connected with connected boundary then

\bigskip
a) $(M\# N)_{-1}=M_{-1}\amalg N_{-1}$.

\bigskip
b) $\partial(Q\amalg Q')=\partial Q\#\partial Q'$

\bigskip
c) $(Q\amalg Q')\times D^1 =(Q \times D^1)\amalg (Q'\times D^1)$

\bigskip
d) $\mathcal{D} (Q\amalg Q')= \mathcal{D}Q\# \mathcal{D}Q'$.

\bigskip
e) $M \# Q= M_{-1}\amalg Q$.

\bigskip
Only part e) seems to require an argument: Let $S$ be the sphere where $M_{-1}$ and $Q_{-1}$ are glued to form $M\# Q$ and $S\times D^1$ a tubular neighborhood of $S=S\times \{0\}$ in $M\# Q$. Take a small disk $D$ in $M\# Q$ around a point in $S\times \{1/2\}$ that intersects $S\times \{1/2\}$ in a disk. Then join $D$ to $\partial Q$ by a thin tube $T$ that does not touch $S\times \{1/2\}$. Then the manifold $R$ which is obtained from $M\# Q$ by removing the interior of $D\cup T$ is clearly diffeomorphic to $M\# Q$. On the other hand, cutting $R$ along the part of $S\times \{1/2\}$ not in $D$ one obtains on one side $M_{-1}$ and on the other a manifold diffeomorphic to $Q$. Joining them back one gets that $R$ is also the union of $M_{-1}$ and $Q$ joined along a disk in each of its boundaries, i.e., $M_{-1}\amalg Q$ and e) is proved.

\bigskip

More surprising are the following facts:

\bigskip

\textbf{Lemma 1 (Fico's Lemmata, [GA]):}\footnote{Unfortunately the Leibniz-Fico
rule: $\mathcal{G}(M \times N)= (\mathcal{G}M \times N)\#(M\times
\mathcal{G}N)$ generalizing 2) is only valid when both $M$ and $N$
are spheres.}

\bigskip

1) \emph{If $M$ and $N$ are connected, then $\mathcal{G}(M\#N)= \mathcal{G}M\#\mathcal{G}N.$}

\bigskip
2) \emph{$\mathcal{G}(S^p \times S^q)=(S^{p+1} \times S^q)\#(S^p \times S^{q+1})$.}

\bigskip

\textbf{Proof: }For 1), just observe from a) and c) above that $$(M\#N)_{-1}\times D^2 = (M_{-1}\amalg N_{-1})\times D^2=(M_{-1}\times D^2)\amalg (N_{-1}\times D^2)$$

\noindent and take the boundary using b).

\bigskip

As for 2), observe that $$(S^p \times S^q)_{-1}=(S^p \times D^q) \cup_{D^p \times D^q }(D^p \times S^q)$$

\noindent (the \emph{plumbing} of two trivial disk bundles). This is not the connected sum $(S^p \times D^q) \amalg (D^p \times S^q)$,
since its boundary is a sphere, the underlying fact being that the two core spheres cannot be separated. However, it is easy to see, by separating the spheres
in the product with $D^1$, that

$$((S^p \times D^q) \cup_{D^p \times D^q }(D^p \times S^q))\times D^1=(S^p \times D^{q+1}) \amalg (D^{p+1} \times S^q)$$

Taking again the product with $D^1$, property c) and the boundary, we obtain 2).

\bigskip

Extending these ideas we have

\bigskip

\textbf{Lemma 2.} \emph{Let $Q,Q'$ be connected $d$-manifolds. Then}

(A)  \emph{If $Q,Q'$ have non-empty boundary then}

$$(Q\# Q')\times D^1=(Q \times D^1) \amalg (Q' \times D^1)  \amalg (D^2 \times S ^{d-1}).$$

$$\mathcal{D}(Q\# Q')=\mathcal{D}Q \# \mathcal{D}Q' \# (S^1 \times S^{d-1})$$

$$\mathcal{TOB}(Q\# Q')=\mathcal{TOB}(Q) \# \mathcal{TOB}(Q') \# (S^2 \times S^{d-1})$$

(B)  \emph{If $Q$ is closed but $Q'$ has non-empty boundary then}

$$\mathcal{D}(Q\# Q')=Q \# (-Q) \# \mathcal{D}Q'$$

$$\mathcal{TOB}(Q\# Q')=\mathcal{G}Q \# \mathcal{TOB}(Q')$$

\textbf{Proof: }first observe that if $Q$ and $Q'$ are disks, $(Q\# Q')$ is diffeomorphic to $D^1 \times S ^{d-1}$. Then, by the argument for property e) preceding Lemma 1, $Q\# Q'$ can be obtained by identifying disks in the boundaries of $Q,Q'$ with corresponding disks in the boundary of $D^1 \times S ^{d-1}$. This \emph{is not} the connected sum along the boundary of those three manifolds because not all their boundaries are connected; actually $Q$ and $Q'$ are attached along different components of the boundary of $D^1 \times S ^{d-1}$. However, this problem
again disappears after multiplying by $D^1$: $Q \times D^1$, $Q' \times D^1$ and the intermediate manifold $D^2 \times S ^{d-1}$ all have connected boundaries
and we have the connected sum along the boundary claimed in A).

Taking its boundary and that of its product with $D^1$ gives the rest of A).

As for B) just observe by e) above that $Q\#Q'= Q_{-1} \amalg Q'$ and take products and boundaries.

\bigskip

\subsection {Cutting vertices.}

\bigskip

Let $P$ be a simple polytope of dimension $d$ with $m$ facets and
$P_\textsf{v}$ be a simple polytope of dimension $d$ obtained from
$P$ by cutting off one vertex $v$.

Let $Z_\textsf{v} =Z(P_\textsf{v})$ be the corresponding connected manifold\footnote{This manifold, the manifold $Z_\textsf{e}$ and any other one obtained by cutting a face of $Z$ can be expressed as an intersection of quadrics by intersecting $P$ with a half space that contains all the vertices of $P$ that do not belong to the face, as explained in section 0.2.}

We will show that the operation $Z_\textsf{v}$ preserves connected sums of sphere products.
However, the result is not simply connected since $P_\textsf{v}$ is never dual $2$-neighborly
because the new created facet does not intersect any of the facets not adjacent to $v$.
So it cannot be combined with the theorems of section 1.
We can partially avoid this problem by the next construction:

Let $Z_\textsf{v}\,'$ be the manifold obtained by duplicating the coefficients of the new variable $x_0$
introduced in the construction of $Z_\textsf{v}$. Notice that this notation is consistent with that of section 1.

\bigskip

\textbf{Theorem 2.1 :}

\bigskip

\textit{1.- $Z_\textsf{v}$ is diffeomorphic to $Z\# Z \# (2^{m-d}-1) (S^1\times S^{d-1}).$}

\bigskip

\textit{2.- $Z_\textsf{v}\,'$ is diffeomorphic to $\mathcal{G}Z \# (2^{m-d}-1) (S^2\times S^{d-1})$.}

\textit{So, if Z is a connected sum of sphere products, any manifold
obtained from it by repeatedly applying the operations $Z_\textsf{v}$
and $Z_\textsf{v}\, '$ is also a connected sum of sphere products.}

\bigskip

\textbf{Proof}: We look at $Z$ as the quotient of $P\times \Z_2^m$ by making
identifications on the facets of $P$. Assume that the vertex $v$
lies in facets $F_1,\dots,F_d$ of $P$ and the rest of the facets are
$F_{d+1},\dots,F_m$. Call $F_0$ the new facet of $P_\textsf{v}$ left
by the removal of a small simplex $\Delta$ around $v$.

Then $Z_\textsf{v}$ is obtained from $P_\textsf{v}\times \Z_2^{m+1}$ by making the
identifications on all of its facets. Let $Z_\textsf{v+}$ be the manifold with boundary
obtained by making all identifications on the facets of $P_\textsf{v}\times \Z_2^m$, except for $F_0$.
Then, in the notation of section 2.1:

$$Z_\textsf{v}=\mathcal{D} Z_\textsf{v+}$$

$$Z_\textsf{v}\, ' =\mathcal{TOB} (Z_\textsf{v+})$$

Now,  $Z_\textsf{v+}$ is $Z$ minus $U$ where $U$ is
obtained from $\Delta \times \Z_2^m$ by making the identifications
on the subsets $\Delta \cap F_i$, where $F_i$ are the facets of $P$.

\bigskip

If we make the identifications on $\Delta \times \Z_2^d$ on the
facets $\Delta \cap F_i$ for $i=1,\dots,d$ we obtain a
disk\footnote{ If we had made identification on $\Delta \times
\Z_2^{d+1}$ on all the facets we would have obtained the sphere
$S^d$.} of dimension $d$, $D^d$. On the product $D^d \times
\Z_2^{m-d}$ we still have to make identifications on the subsets
$\Delta \cap F_i$ for $i=d+1 \dots m$, but these
intersections are empty. Therefore $U=D^d \times \Z_2^{m-d}$ and
$Z_\textsf{v+}$ is $Z$ with $2^{m-d}$ disks removed.

\bigskip

Now we use the lemma:

\textbf{Lemma.} \textit{Let $M^d$ be a connected manifold and
$M_{-k}$ be obtained from $M$ by removing $k$ open disks. Then}

\bigskip

1) $M_{-k} \times D^1= (M_{-1}\times D^1) \scf (k-1)(D^2\times S^{d-1}).$

\bigskip

2) $\mathcal{D} (M_{-k}) = M\# (-M) \# (k-1)(S^1\times S^{d-1}).$

\bigskip

3) $\mathcal{TOB} (M_{-k}) = \mathcal{G}M \# (k-1)(S^2\times S^{d-1}).$

\bigskip

\textbf{Proof:} Since $M$ is connected, we can assume all the open disks lie in
the interior of a closed disk and therefore $M_{-k}$ can be seen as
$M_{-1}$ with $k-1$ copies of $D^1\times S^{d-1}$ attached along disks in its boundary.
This is \emph{not} a connected sum along the boundary, since the boundary of $D^1\times S^{d-1}$ is not
connected, but after multiplying times $D^1$ this problem disappears
and $M_{-k} \times D^1$ is the connected sum claimed in 1). Taking its boundary
we get part 2) and taking its boundary after multiplying again by $D^1$ gives
part 3). So the lemma is proved and so are parts 1) and 2) Theorem 2.1.

The last part is obvious for the operation $Z_\textsf{v}$ and it follows from
Fico's Lemmata (section 2.1) for the operation $Z_\textsf{v}\, '$.

\bigskip

\textbf{Remarks:}

\bigskip

Although the combinatorial type of $P_\textsf{v}$
depends on which vertex is cut off, the diffeomorphism type of
$Z_\textsf{v}$ does not. This gives many examples of
different polytopes that produce the same manifold $Z$.

The analogous facts about $Z^\C$ appear in [B-M].

\bigskip

In particular, if $P$ is obtained from a simplex by succesively cutting off vertices
then $Z(P)$ is a connected sum of sphere products. Such a polytope is a
\emph{dual stack polytope}, and its dual can be obtained from a simplex by
succesively building up pyramids on facets.

\bigskip

Its diffeomorphism type depends only on the number of times each
operation is applied and not on which vertices, or in which order, they are applied.
As an exercise, write the formula for the effect on the manifold $Z$
of successively cutting off $n$ vertices from $P$.

\bigskip

\bigskip

Part 2 of Theorem 2.1 can be combined with the theorems of section 1 to obtain
that connected sums are preserved when applying repeatedly the $Z^J$ and the $Z_\textsf{v}\, '$ constructions. A more complete version of this will appear as Theorem 2.4 below.

\bigskip

Nevertheless, one would like to have a much more general vertex
cutting result. We do not know, for example, if
$Z^\C(P_\textsf{v}\,)$ is a connected sum of sphere products
whenever $Z^\C(P)$ is.

More generally, we can make the following

\bigskip

\textbf{Conjecture:}

\bigskip

\textit{$Z^\C(P_\textsf{v}\,)$ is diffeomorphic to $\mathcal{G}
(Z^\C(P)) \# \sc_{j=1}^{m-d}{m-d \choose j} (S^{j+2}\times S^{d+m-j-1})$.}

\bigskip

We now prove it under a restriction on the number of facets of the polytope:

\bigskip

\textbf{Theorem 2.2:} \textit{If $P$ is a polytope of dimension $d$ with $m < 3d$ facets then the above conjecture is true.}

\bigskip

We will follow the same strategy as in the previous results, although now the steps are less elementary:

\bigskip

1) $Z^\C(P)$ is obtained from $P \times (S^1)^m$ by making identifications on the $m$ facets. To obtain $Z^\C(P_\textsf{v})$
we have to multiply $P_\textsf{v}$ by $(S^1)^{m+1}$ and make the identifications on the $m+1$ facets. Let $Z^\C_+(P_\textsf{v})$
be the manifold with boundary obtained from the product $(P_\textsf{v})\times (S^1)^m$ by making the identifications on all
the facets except the new facet $0$. Then $Z^\C(P_\textsf{v})$ is obtained from $Z^\C_+(P_\textsf{v})$ by multiplying by
$S^1$ and making the identifications on its boundary. In the notation of section 2.1 this means

$$Z^\C(P_\textsf{v})=\mathcal{TOB}(Z^\C_+(P_\textsf{v})).$$

So we have to understand $Z^\C_+(P_\textsf{v})\times D^1$. Now $Z^\C_+(P_\textsf{v})$ is obtained from $Z^\C(P)$
by removing the manifold $U$ obtained from $\Delta^d\times (S^1)^m$ by making the same identifications as for $Z^\C(P)$, where
$\Delta^d$ is the simplex removed to create the new facet $0$. If we make the identifications on $\Delta^d\times (S^1)^d$
corresponding to the facets of $\Delta^d$ that are part of facets of $P$ we obtain the disk $D^{2d}$. On the product $\Delta^d\times (S^1)^m$ we still have to make
identifications on the intersections with the $m-d$ facets not adjacent to the vertex $v$ but these intersections are empty.
Therefore $U=D^d \times (S^1)^{m-d}$ and $Z^\C_+(P_\textsf{v})$ is $Z^\C(P)$ with a thickened $(m-d)$-torus removed.

\bigskip

2) This $(m-d)$-torus can be contracted to a point in $Z^\C(P)$ because it is an orbit of the action of $T^m$ (see section 0.2).
Since $2(m-d)< m+d $ it is isotopic to an $(m-d)$-torus inside an open disk in $Z^\C(P)$. Therefore
$$Z^\C_+(P_\textsf{v})=Z^\C(P) \# (S^{m+d} \setminus (T^{m-d}\times int(D^{2d}))).$$

\bigskip

Again, since $2(m-d)< m+d $ we can apply Theorem A2.3 of the Appendix to the last summand and Lemma 1 of section 2.1 to the whole sum to obtain the Theorem.

\bigskip

Without the restriction on the number of facets of $P$, it can be proved that the $(m-d)$-torus can be engulfed by an open disk in $Z^\C(P)$.
To prove the general conjecture, one would need to prove that this torus is in some sense \emph{standard} inside that disk and some version of Theorem A2.3 of Appendix A3 for that type of \emph{standard} torus.

A more general formula can be conjectured when the construction $Z^J$ is immediately applied after cutting off a vertex. Again a partial case can be deduced from Theorem A2.2 of the Appendix.

\bigskip 

We need also, of course, a way to include all the non-simply
connected manifolds into the picture.

\bigskip

\textbf{Example: the truncated cube.}

\bigskip

We will study now the cube $C=I\times I\times I$ and the truncated cube $C_\textsf{v}$.

\begin{center}
\includegraphics[width=4in]{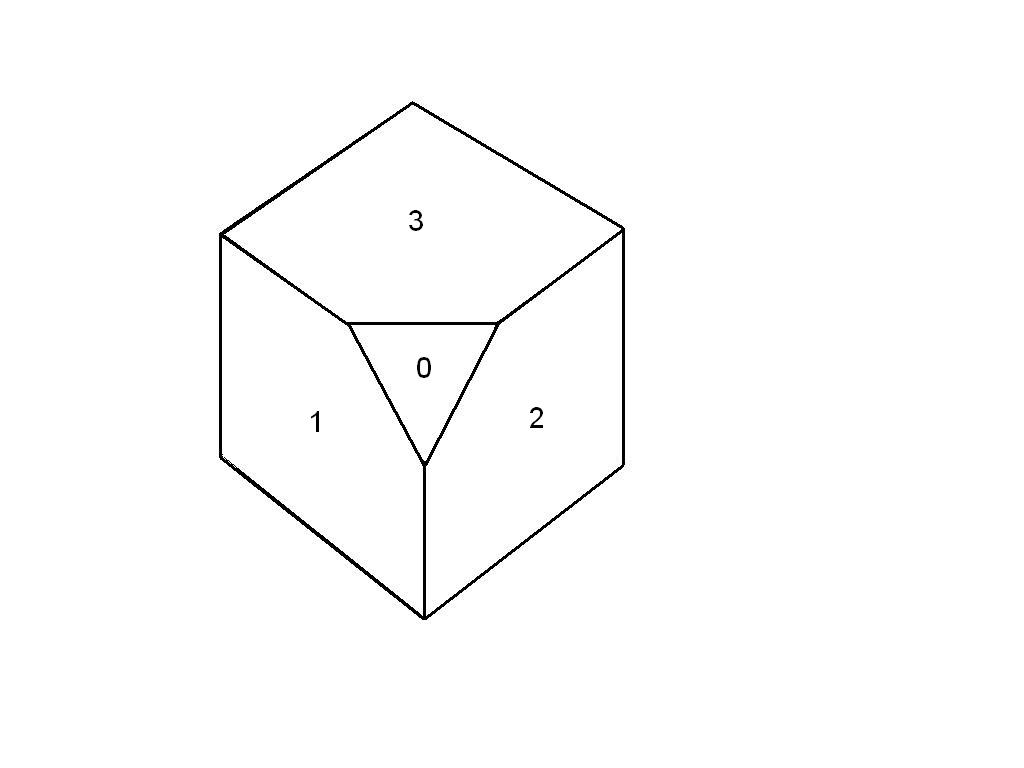}
\end{center}

Now we have $$Z(C)=S^1 \times S^1 \times S^1 $$

$$Z^\C(C)=S^3 \times S^3 \times S^3 $$

From Theorem 2.1 we know that

$$Z(C_\textsf{v})=(S^1 \times S^1 \times S^1) \# (S^1 \times S^1 \times S^1) \# 7 (S^1 \times S^2)$$

From Theorem 2.2  we know that

$$Z^\C(C_\textsf{v})= \mathcal{G} (S^3 \times S^3 \times S^3) \# 3 (S^3 \times S^7) \# 3 (S^4 \times S^6) \# (S^5 \times S^5).$$

It follows from the cohomology ring of $\mathcal{G} (S^3 \times S^3 \times S^3)$ that this manifold does not decompose into a non-trivial connected sum: by the argument of the Proposition in section 3, the 3-dimensional generators would have to lie in one of the summands and so would the rest of the generators of dimension less that 10 since they are related by non-trivial products with the first ones, so one of the summands must be a homotopy sphere. This shows that we have reached the whole decomposition of $Z^\C(C_\textsf{v})$ as a connected sum.

\bigskip

This answers another question in [B-M], p.111. It will prove useful in section 3.

\bigskip

We can also describe the topology of the manifolds associated with the truncated $d$-cube. In fact, Theorem 2.2 gives us the topology of $Z^\C$ of an iterated truncation of up to $d$ vertices
of the $d$-cube.

\subsection {Cutting edges and other faces.}

Let $P$ be a simple polytope of dimension $d$ and $m$ facets and
$P_\textsf{e}$ be a simple polytope of dimension $d$ obtained from
$P$ by cutting off one edge $e$.

\bigskip

Let $Z_\textsf{e}=Z(P_\textsf{e})$. Again, the resulting manifold is not simply connected, but once more we can prove that it becomes simply connected after applying to it the construction $Z'$ based on the new facet of $P_\textsf{e}$. Let $Z_\textsf{e}\, '$ be that manifold.

\bigskip

\textbf{Theorem 2.3:} \textit{If $Z$ is 1-connected (or, equivalently, if $P$ is dual 2-neighborly) then}

1) \textit{$Z_\textsf{e}$ is diffeomorphic to $$Z\# Z \# (2^{m-d-1})(S^2\times S^{d-2})\# (2^{m-d-1}-1)(S^1\times S^{d-1}).$$}

2) \textit{$Z_\textsf{e}\, '$ is diffeomorphic to $$\mathcal{G}(Z) \# (2^{m-d-1})(S^3\times S^{d-2})\# (2^{m-d-1}-1)(S^2\times S^{d-1}).$$}

\bigskip

\textbf{Proof:} We look at $Z$ as the quotient of $P\times \Z_2^m$ by making
identifications on the facets of $P$. Denote by $F_1,\dots,F_m$ these facets, where $F_1,\dots F_{d+1}$ are the ones touched by $e$. Call $F_0$ the new facet of $P_\textsf{e}$ left by the
removal of a regular neighborhood $e\times \Delta^{d-1}$ of $e$.

Then $Z_\textsf{e}$ is obtained from $P_\textsf{e}\times \Z_2^{m+1}$
by making the identifications on its $m+1$ facets.
Let $Z_\textsf{e+}$ be obtained from $P_\textsf{e}\times \Z_2^m$, by making identifications on the facets
$F_1,\dots, F_m$. Then clearly we have

$$Z_\textsf{e}= \mathcal{D}(Z_\textsf{e+}) $$

Also,  $Z_\textsf{e}\, '{}$ is obtained from $P_\textsf{e}\times \Z_2^m \times S^1$
by making the identifications on its facets, so, in the notation of section 2.1,

$$Z_\textsf{e}\, '= \mathcal{TOB}(Z_\textsf{e+}) $$

So we have to describe $Z_\textsf{e+}\times D^1$ and $Z_\textsf{e+}\times D^2$ in order to know their boundaries.

Now $Z_\textsf{e+}$ is $Z$ minus $U$ where $U$ is obtained from
$e\times \Delta^{d-1} \times \Z_2^m$ by making the identifications on its intersections with
the facets $F_1,\dots, F_m$ of $P$.

If we make first the identifications on $e\times \Delta^{d-1}\times \Z_2^{d+1}$ on all its facets we obtain the product $S^1 \times S^{d-1}$,
but if we only do it on those facets that are part of the facets of $P$ (i.e., excluding the facet $0$) we obtain
$S^1 \times D^{d-1}$. On the product $e\times \Delta^{d-1}\times \Z_2^m$
we still have to make identifications on the intersections with the facets
$i$ for $i=d+2, \dots, m$ but these
intersections are empty. Therefore $U=S^1 \times D^{d-1} \times \Z_2^{m-d-1}$ and
$Z_\textsf{e+}$ is $Z$ with $2^{m-d-1}$ copies of $S^1 \times D^{d-1}$ removed.

Then we have the lemma\footnote{Recall the notation of section 2.1: $M_{-1}$ denotes $M$ minus an open disk. For dimension 3 the lemma is not true, but the theorem is, since the only dual 2-neighborly polytope is the simplex.}:

\bigskip

\textbf{Lemma.} \emph{Let $M^d$ be a simply connected manifold of dimension at least $4$ and $M_{\sim k}$ be obtained from $M$ by removing $k>0$ copies of $S^1 \times D ^{d-1}$.}

\emph{Then $M_{\sim k} \times D^1$ is diffeomorphic to} $$M_{-1}\times D^1  \amalg k(D^3\times S^{d-2}) \amalg (k-1)(D^2\times S^{d-1}).$$

\bigskip

 \textbf{Proof:} Since $M$ is simply connected, we can assume all copies of $S^1 \times D^{d-1}$ lie inside an
open disk, so we have only to prove it for $M=S^d$.

In this case the result is true for $k=1$ (even before multiplying by $D^1$) since $S^d$
minus a copy of $S^1 \times D^{d-1}$ is diffeomorphic to $D^2 \times S^{d-2}$.

Inductively, $M_{\sim (k+1)}$ is obtained from $M_{\sim k}$ by removing from it one copy of $S^1 \times D^{d-1}$ so $M_{\sim (k+1)}=M_{\sim k} \# M_{\sim 1}$.

Now Lemma 2 A) of section 2.1 gives the induction step:

$$M_{\sim (k+1)}\times D^1=(M_{\sim k}\# M_{\sim 1})\times D^1=(M_{\sim k} \times D^1)\amalg (M_{\sim 1} \times D^1)\amalg (D^2 \times S^{d-1})=$$
$$(M_{-1}\times D^1) \amalg k(D^3\times S^{d-2})\amalg (k-1)(D^2\times S^d)\amalg (D^3 \times S^{d-2})\amalg (D^2 \times S^{d-1})$$
$$=(M_{-1}\times D^1) \amalg (k+1)(D^3\times S^{d-2})\amalg k(D^2\times S^d).$$

\noindent and the Lemma is proved.

\bigskip

Applying it to $Z$ we obtain that $Z_\textsf{e+}\times D^1$ is diffeomorphic to $$(Z_{-1}\times D^1) \amalg 2^{m-d-1}(D^3\times S^{d-2})\#(2^{m-d-1}-1)(D^2\times S^{d-1}).$$

\bigskip

Taking its boundary we obtain 1). Taking the boundary of its product with $D^1$ we obtain 2).

\textbf{Remarks:}

1.- Although the combinatorial type of $P_\textsf{e}$ depends on which edge is cut off, the diffeomorphism type of $Z_\textsf{e}$
does not. This gives many new examples of different polytopes that produce the same manifold $Z$.

2.- If $Z$ is not simply connected then Theorem 2.3 and the previous observation are not true. For example, if $P$ is
the triangular prism $\Delta^2\times I$, then by cutting off an edge we can obtain either the pentagonal book or the cube, whose
corresponding manifolds are $\# 5(S^1\times S^2)$ and $S^1\times S^1\times S^1$, respectively. The reason lies in the fact that when the edge connects two facets of $P$ that do not intersect, then this edge produces an embedded $S^1$ that is not trivial in $H_1(Z)$ and therefore cannot be contained inside a disk. It remains to be clarified what happens in this case.

3.- One problem with part 1 is that it cannot be immediately iterated since the resulting manifold $Z_\textsf{e}$ is not simply connected.

\bigskip

Combining Theorem 2.3 with our previous theorems we get

\bigskip

\textbf{Theorem 2.4:} \textit{If $Z = Z(P)$ is a connected sum of sphere products and is
simply connected of dimension at least 5, then any manifold
obtained from $Z$  by repeatedly applying the $Z^J$, $Z_\textsf{v}\, '$ and $Z_\textsf{e}\, '$ constructions
(any number of times each and in any order) is also a connected sum of sphere products.\footnote{This includes, of course, the corresponding moment-angle manifolds.}}

\bigskip

The basic manifold can be any sphere, any $Z(P)$ with free homology for a dual-neighborly polytope (Theorems 1.3 and 1.4) or any simply connected manifold with $m - d = 3$ of dimension at least 5 that is not a triple sphere product ([LdM1]).

It is known that there are many dual neighborly polytopes in each dimension. If to these we add the infinite family stemming from each of them by applying the three constructions above (which produce polytopes which are \emph{not} dual neighborly) we can conclude that there is a large number of polytopes $P$ for which $Z(P)$ is a connected sum of sphere products.

\bigskip

The theorems on this section can be generalized to the process of removing any simplicial face of $P$ but need, of course, stronger connectivity hypotheses.

One can conceive also results about removing non-simplicial faces of $P$. It has been conjectured that any simple polytope can be realized as a face of a dual neighborly polytope, so in principle we need to consider the effect of cutting off any possible simple polytope.

\newpage

\section {On the cohomology ring of $Z(P)$.}

In [B-M] the cohomology ring of $Z^\C(P)$ is described (by an explicit rule) in terms of intersection products on the relative homology groups of pairs
$(P,P_J)$ where $P_J$ is a union of facets of $P$. A slightly different rule, valid for general moment-angle complexes, appears in [Ba]. We shall see now
by examples that these rules must be drastically modified in the case of the manifolds $Z(P)$.

We will show that the cohomology rings of the manifolds associated to the truncated cube are not isomorphic as \emph{ungraded rings}.

\bigskip

\emph{We will consider the cohomology rings with $\Z_2$ coefficients.} Since all the homology groups of the spaces involved are free,
this will be enough to establish the difference without having to compute signs, so in fact we have a stronger statement.

\bigskip

We have shown at the end of section 2.2 that

$$Z(C_\textsf{v})=(S^1 \times S^1 \times S^1) \# (S^1 \times S^1 \times S^1) \# 7 (S^1 \times S^2)$$

$$Z^\C(C_\textsf{v})= \mathcal{G} (S^3 \times S^3 \times S^3) \# 3 (S^3 \times S^7) \# 3 (S^4 \times S^5) \# (S^5 \times S^5).$$

We begin by describing the cohomology ring of the first summands of each manifold:

$$X:= (S^1 \times S^1 \times S^1) \# (S^1 \times S^1 \times S^1)$$
$$Y:=\mathcal{G} (S^3 \times S^3 \times S^3)$$

$H^*(X)$ has:

\bigskip
One generator $1$ in dimension 0.

\bigskip
6 generators $a_1, a_2, a_3, a_1', a_2', a_3'$ in dimension 1.

\bigskip
6 generators $b_{1,2},b_{1,3}, b_{2,3}, b_{1,2}',b_{1,3}', b_{2,3}'$ in dimension 2.

\bigskip
One generator $f$ (the fundamental class) in dimension 3.

\bigskip
The multiplication rules are:

\begin{center}
$a_i a_j = b_{i,j}$ for $i\ne j$

\bigskip
$a_1 a_2 a_3 = f$

\bigskip
$a_i' a_j' = b_{i,j}'$ for $i\ne j$

\bigskip
$a_1' a_2' a_3' = f$
\end{center}

\noindent all other products not derived from these being $0$.

\bigskip

On the other hand $H^*(Y)$ has:

\bigskip
One generator $1$ in dimension 0.

\bigskip
3 generators $A_1, A_2, A_3$ in dimension 3.

\bigskip
3 generators $A_1', A_2', A_3'$ in dimension 4.

\bigskip
3 generators $B_{1,2},B_{1,3}, B_{2,3}$ in dimension 6.

\bigskip
3 generators $B_{1,2}',B_{1,3}', B_{2,3}'$ in dimension 7.

\bigskip
One generator $F'$ (the fundamental class) in dimension 10.

\bigskip
(In this case a generator is denoted by $C'$ when it is inherited from the element $C\times u$ in the cohomology
of $(S^3 \times S^3 \times S^3)\times S^1$ where $u$ is the generator of $H^1(S^1)$).

\bigskip
The multiplication rules are:

\begin{center}
$A_i A_j = B_{i,j}$ for $i\ne j$

\bigskip
$A_i A_j' = B_{i,j}'$ for $i\ne j$

\bigskip
$A_i A_j A_k' = F'$ for $i,j,k$ all different.

\end{center}

\noindent all other products not derived from these being $0$.

\bigskip

Observe in particular that

$$ A_i' A_j' =0\,\,\, (*)$$

\noindent for all $i,j$, simply because the product lies in $H^8(Y)=0$.

This is the main difference with the previous ring: three indecomposable, linearly independent elements
have all their products equal to zero. To make this idea into a formal proof of non-isomorphism we have to work modulo the
decomposable elements:

In both rings there is the maximal ideal $\mathcal{M}$ of non-units, spanned by all the generators except 1.
The square of this ideal $\mathcal{M}^2$ is generated by the decomposable elements and its cube is generated
by the fundamental class. There is an induced multiplication

$$M:\mathcal{M}/\mathcal{M}^2\times \mathcal{M}/\mathcal{M}^2 \rightarrow \mathcal{M}^2/\mathcal{M}^3$$

An ungraded ring isomorphism between the two rings would have to preserve the above ideals and multiplication.

\bigskip

\textbf{Proposition.} \emph{For the ring $H^*(X)$, any three elements in $\mathcal{M}/\mathcal{M}^2$ all whose products are zero
in $\mathcal{M}^2/\mathcal{M}^3$ must be linearly dependent.}

\bigskip

\textbf{Proof:} Let $$x_j=\mathop{\Sigma}\limits_{i=1}^3 \lambda_j^i a_i+\mathop{\Sigma}\limits_{i=1}^3 \mu_j^i a_i'$$

\noindent for $i=1,\dots, 3$ be (representatives of) such elements. Then

$$x_j x_k= \mathop{\Sigma}\limits_{r<s} (\lambda_j^r  \lambda_k^s + \lambda_j^s  \lambda_k^r) b_{r,s} + \mathop{\Sigma}\limits_{r<s} (\mu_j^r  \mu_k^s + \mu_j^s  \mu_k^r) b_{r,s}'.$$

The fact that this product is 0 means that all $2\times 2$ determinants of the matrices $\lambda_j^i$ and $\mu_j^i$ are 0. This means that in each there is at most one linearly independent column and that the matrix of coefficients of the $x_i$,  which is formed by the columns of those two, has rank at most 2. Thus proposition 4.1 is proved.

\bigskip

To complete the rings of $Z(C_\textsf{v})$ and $Z^\C(C_\textsf{v})$ one has to add in each case to the corresponding previous one 7 couples of elements $u_i,v_i$ such that each
product $u_i v_i$ is the fundamental class and the product of each $u_i,v_i$ with all the rest of generators is 0. We have now:

\bigskip

\textbf{Theorem 3.1} \emph{The rings $H^*(Z(C_\textsf{v}),\Z_2)$ and $H^*(Z^\C(C_\textsf{v}),\Z_2)$ are not isomorphic as ungraded rings.}

\bigskip

\textbf{Proof:} We follow the same lines as for the previous proposition: now the rings considered there must be complemented with 7 couples of generators $u_i,v_i$ such that each
product $u_i v_i$ is the fundamental class and the product of each $u_i$ and $v_i$ with all the rest of generators is 0. Since the fundamental class is in $\mathcal{M}^3$ this means that  $H^*(Z^\C(C_\textsf{v}))$ has 17 generators in $\mathcal{M}/\mathcal{M}^2$ all whose products are 0 under the induced multiplication $M$. If in the case of $H^*(Z(C_\textsf{v}),\Z_2)$ we have 17 elements in $\mathcal{M}/\mathcal{M}^2$ for which all products under $M$ are 0, we can express them in terms of the 20 generators $a_i,a_i',u_i,v_i$ by a $17\times20$ matrix. By the same computation as above, at most two of the first six columns of this matrix are independent, so the matrix has at most rank 16 and the 17 elements are not linearly independent. This proves the Theorem.

\bigskip

\textbf{Remarks.} 1.- Thus the analogues of the product rules in [Ba], [B-M] for the ring $H^*(Z^\C(P))$ are not valid in general for the ring $H^*(Z(P))$.

\bigskip

2.- Theorem 3.1 contradicts the results on the cohomology ring of the complement of a real coordinate subspace arrangements in [dL] (see Remark 4.1) since everyone of our manifolds has the homotopy type of such a complement. Actually, the example where $P$ is a pentagon (the \textit{truncated square}) already shows that the product rule has to be modified in the real case. But for the pentagon the cohomology rings over $\Z_2$ turn out to be isomorphic as ungraded rings, even if the \textit{natural} additive isomorphism does not preserve the ring structure. But in our example there is no possible ring isomorphism at all, and the same happens for any truncated cube of dimension greater that 1. We thank Taras Panov for pointing out that another example can be found in [G-P-W]. This example is a subspace arrangement that does not correspond to a convex polytope.

\bigskip

3.- It can be shown that the cohomology rings of $Z^\C(P)$ and $Z^\C(P')$ are isomorphic as \emph{ungraded rings}. (See [B-B-C-G3]); this should also be a consequence of the other known descriptions of the cohomology rings ([B-M], [B-P]). For example, if $P$ is a product of two simplices then $Z^\C(P) = S^{2p-1}\times S^{2q-1}$ and $Z^\C(P') = S^{2p+1}\times S^{2q-1}$ and the ungraded cohomology rings are isomorphic, with generators $\{1,x,y,xy\}$. Or, if $Z^\C(P)$ is a connected sum of sphere products then $Z^\C(P')$ is also one, with the same number of summands, but where inside each summand one of the factors has risen its dimension by two.

\bigskip

The previous Theorem implies that this also cannot be true in general for the cohomology rings of $Z(P)$ and $Z(P')$, not even when taken with $\Z_2$ coefficients, since the passage from $Z(P)$ to $Z^\C(P)$ goes through the iteration of the operation $Z(P')$.

We claim that the correct rules are as follows: Let $L$ be a collection of facets of $P$, we denote by $P_L$ the union of the facets of $P$ in $L$.

\bigskip

\textit{Intersection product: There is an isomorphism }$$H_i(Z) \approx \mathop{\oplus}\limits_L H_i(P,P_L)$$

\noindent \textit{so that the intersection of two homology classes in $Z$ is given by the sum of the intersection products on the right hand side}

$$H_i(P,P_L) \otimes H_j(P,P_J) \rightarrow H_{i+j-d}(P,P_{L\cap J}).$$

\bigskip

\textit{Cup product: There is an isomorphism }$$H^i(Z) \approx \mathop{\oplus}\limits_L H^i(P,P_L)$$

\noindent \textit{so that the cup product of two cohomology classes in $Z$ is given by the sum of the cup products on the right hand side}

$$H^i(P,P_L) \otimes H^j(P,P_J) \rightarrow H^{i+j}(P,P_{L\cup J}).$$

In other words, the rule that a product must be 0 unless $L$ and $J$ satisfy a strict set-theoretical condition, must be dropped.

\bigskip

Details and extensions will appear in a forthcoming article ([G-LdM]).

\newpage

\section*{Appendix.}

\subsection* {A1. Recognizing a connected sum of sphere products.}

Recall that, after Thom, Smale and Milnor, the best way to identify a compact manifold $M$ is to study an adequate compact manifold $Q$ whose boundary
is $M$. For example, if $M^n$ is a homotopy sphere, we cannot tell if it is diffeomorphic to the standard sphere $S^n$ by looking
at the usual topological invariants of $M$: since these are homotopy invariants there is nothing to look at. But we know that $M$ is the boundary of some manifold
$Q$. The first thing to ask is if $Q$ is parallelizable. If there is no such parallelizable $Q$ then $M$ is not standard since $S^n$ bounds a disk. If one finds
a $Q$ that is parallelizable, then one can decide if $M$ is standard by studying topological invariants of $Q$ ([K-M]).

Now a connected sum of sphere products

$$\sc^k_{i=1} (S^{p_i}\times S^{d-p_i}) $$

\noindent is the boundary of the connected sum along the boundary

$$\scf_{i=1}^k (S^{p_i}\times D^{d-p_i+1}) $$

\noindent which has the following properties: It is simply connected if $1<p_i$, has simply connected boundary if $1<p_i<d-2$, its homology
groups are free and bases of them can be represented by disjoint embedded spheres with trivial normal bundle.

\bigskip

To determine such a connected sum it is enough to verify those conditions:

\bigskip

\textbf{Theorem A1.} \emph{Let $Q$ be a manifold with boundary satisfying:}

\bigskip

a) \emph{$Q$ is of dimension $d+1 \ge 6$}

\bigskip
b) \emph{$Q$ is simply connected with simply connected boundary.}

\bigskip
c) \emph{$H_i(Q)$ is free with basis $\alpha_{i,j}$ for all $i$ and $H_i(Q)=0$ for $i\ge d-1$}

\bigskip
d) \emph{There is a collection $\{S^i_j\}$ of disjoint embedded spheres with trivial normal bundle inside $Q$ that represent the basis elements $\alpha_{i,j}$ of $H_i(Q)$.
}
\bigskip

\emph{Then $Q$ is diffeomorphic to a connected sum along the boundary
}
$$\scf_{i=1}^k (S^{p_i}\times D^{d-p_i+1}) $$

\emph{\noindent and therefore $\partial Q$ is diffeomorphic to }

$$\sc_{i=1}^k (S^{p_i}\times S^{d-p_i}).$$

\textbf{Proof:} For every sphere $S^i_j$ take a closed product neighborhood $U_{i,j}=S^i_j \times D^{d-i}$ in such a way that they are all disjoint.
Connect all $U_{i,j}$ by a minimal collection of thin tubes and call $Q_1$ the union of all the $U_{i,j}$ and their connecting tubes, which we can
assume is contained in the interior of $Q$. Then $Q_1$ is a connected sum of the required form for $Q$ and we will prove that these two manifolds are
diffeomorphic. Let $H$ be equal to $Q$ minus the interior of $Q_1$, so that $H$ is a cobordism between the manifolds $\partial Q$ and $\partial Q_1$.

Now $Q_1$, $\partial Q_1$ are simply connected because all the spheres $S^i_j$ have dimension $2\ge i\ge d-2$ and so is $H$ because $Q$ is simply connected
and every loop in $H$ bounds a disk in $Q$ that can be made disjoint from $Q_1$ since the spheres have codimension greater than 2. Now $H_i(H,\partial Q_1)=H_i(Q, Q_1)=0$
for all $i$ because the inclusion of $Q_1$ in $Q$ induces isomorphisms of homology groups. Therefore $H$ is a simply connected $h$-cobordism of dimension
at least 6 between simply connected manifolds and so is a product by the $h$-cobordism theorem. This proves that $Q$ is diffeomorphic to $Q_1$ and the
theorem is proved.

\subsection*{A2. Complements of sphere products in spheres.}

We will work in the unit sphere $S^n$ in $\R^{n+1}$ and assume we have a decomposition

$$\R^{n+1}=\R^{n_1+1}\times\dots \R^{n_k+1}\times \R^q.$$

We have then $n=\mathop{\Sigma}\limits_{i=1}^k n_i+q+k-1$ and we will denote an element in $\R^{n+1}$ as a $(k+1)$-tuple
$(X_1,\dots ,X_k, Y)$ with $X_i\in \R^{n_i+1}$ and $Y\in \R^q$.

\bigskip

%We will denote by $\{e^i_j\}$ for $j=0,\dots ,n_i$ the standard basis of $\R^{n_i+1}$.

%\bigskip

We will consider the product of spheres $$\mathcal{P}=S^{n_1}\times\dots S^{n_k}\times \{0\}$$

\noindent given by the equations

$$ X_i\cdot X_i = \frac{1}{k},\,\,\,\, Y=0.$$

\bigskip

$\mathcal{P}$ has dimension $\mathop{\Sigma}\limits_{i=1}^k n_i$ and codimension $q+k-1$.

\bigskip

Within each sphere $S^{n_i}$ we consider the base point $p_i=(\frac{1}{\sqrt{k}},0,\dots,0)$. We denote by $p$ the point $(p_1,\dots,p_k,0)\in \R^{n+1}$

%\frac{1}{\sqrt{k}} e^i_0=

\bigskip

We are interested in the complement of an open tubular neighborhood of $\mathcal{P}$ in $S^n$. We will denote this complement by $\mathcal{E}$.

\bigskip

Let $I=\{1,\dots, k\}$. If $J\subset I$ we denote by $\left|J\right|$ the number of elements of $J$ and
by $n_J$ the sum $\mathop{\Sigma}\limits_{i\in J}n_i$.

\bigskip

We will make use of the subproduct $\mathcal{P}_J$ of $\mathcal{P}$ defined by the equations

$$X_i=p_i \,\,\,\, i\notin J$$ $$Y=0.$$

So $\mathcal{P}_J$ is diffeomorphic to the product $\Pi_{i\in J}S^{n_i}$ of dimension $n_J$,
the rest of the coordinates having fixed values.

\bigskip

If $n_i>0$ for $i=1,\dots , k$ then $\mathcal{P}$ is connected and its reduced homology groups are free with a generator of dimension $n_J>0$ represented by the fundamental class of $\mathcal{P}_J$ for every non-empty $J\subset I$.

We will show that if $n_i>0$ for $i=1,\dots , k$ we have Alexander dual generators for the homology of $\mathcal{E}$ represented by embedded round spheres. The Alexander dual of a single sphere $S^{n_i}$ will be a great sphere of dimension $n-n_i-1$ lying in a complementary subspace. The Alexander dual of the products of spheres will be spheres that become smaller and more slanted as the number of factors grows. The Alexander dual of the whole product $\mathcal{P}$ will be a small sphere which is a fiber of its spherical normal bundle. 
\newpage

A schematic representation of the construction is the following:

\begin{center}
\includegraphics[width=4in]{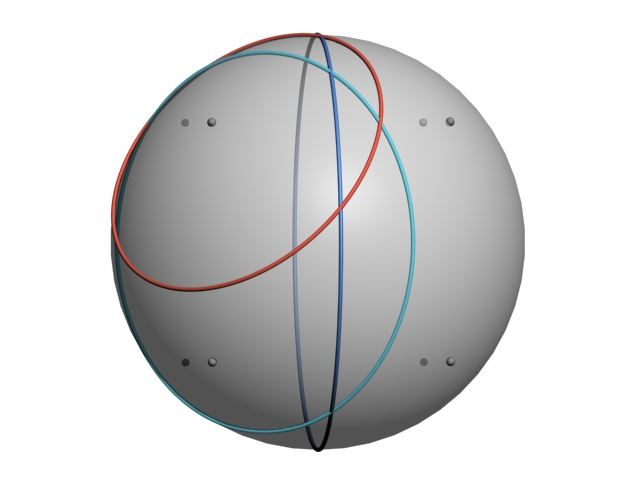}
\end{center}

\vspace{-10pt}

Here we have the picture for $n=2$, $k=3$, $n_1=n_2=0$ (a case \emph{not} covered by the hypothesis of Theorem A2.1 below). The product $\mathcal{P}=S^0\times S^0\times S^0$ is given by the 8 marked points. The vertical great circle is the Alexander dual of an $S^0$ given by one point on each vertical hemisphere. The smaller vertical circle is the Alexander dual of an $S^0\times S^0$ consisting of the points in the upper hemisphere since the small disk this circle bounds contains exactly one of its points. The smallest slanted circle is the Alexander dual of the whole product, since it surrounds just one of its points.

\bigskip

\textbf{Theorem A2.1}. \textit{Assume $n_i>0$ for $i=1,\dots , k$. Then there is a basis for the homology of $\mathcal{E}$ represented by embedded round spheres. This basis is in each dimension the Alexander dual basis of the basis $\mathcal{P}_J$ of the reduced homology of $\mathcal{P}$ in the corresponding dimension.}

\bigskip

\bigskip

\textbf{Proof:} The result is well-known for $k=1$ so we assume $k>1$.

We will consider for each $J\subset I$ the following subspaces of $\R^{n+1}$:

\bigskip

$L_J$ will be the vector subspace of $\R^{n+1}$ of codimension $n_J$ consisting of those points $(X,Y)$ such that, for $i\in J$, $X_i$ is a multiple of $p_i$:

$$X_i=t_i p_i\,\,\,\,\,\,\,\, for\,\, i\in J$$

Inside its unit sphere $L_J\cap S^n$ we consider the set $D_J$ of points satisfying

$$\mathop{\Sigma}\limits_{i\in J}\,t_i\ge \left|J\right|-1.$$

\noindent and the set $S_J$ of points satisfying

$$\mathop{\Sigma}\limits_{i\in J}\,t_i=\left|J\right|-1.$$

Since the function $\mathop{\Sigma}\limits_{i\in J}\,t_i$ takes at $p\in L_J\cap S^n$ the value $\left|J\right|$, it follows that

\bigskip

1) $D_J$ is a disk with non-empty interior in $L_J\cap S^n$.

\bigskip

2) $S_J$ is a round sphere which is the boundary of $D_J$.

%The extreme values this function achieves in $S^n$ are $\pm \sqrt{k \left|J\right|} \ge  \left|J\right|$, taken  at $X=\pm \mathop{\Sigma}\limits_{i\in %J}\, t p_i$ where $t= \sqrt{k/\left|J\right|}$.

\bigskip

We claim that the $S_J$ represent the Alexander duals of the basis $\mathcal{P}_J$ of the homology of $\mathcal{P}$. We show this by proving that for each $J\ne \emptyset$ :

\bigskip

A) $S_J$ is contained in $\mathcal{E}$.

\bigskip

B) The disk $D_J$ in $S^n$ whose boundary is $S_J$ intersects $\mathcal{P}_J$ transversely in exactly the point $p$.

\bigskip

C) For any $J'\ne J$ with $n_J=n_{J'}$ the intersection number of $D_J$ with $\mathcal{P}_{J'}$ is zero.

\bigskip

To prove A) and B) first observe that a point $(X,Y)\in \mathcal{P} \cap (L_J\cap S^n)$ must satisfy $Y=0$ and both $X_i\cdot X_i = \frac{1}{k}$ and $X_i=t_i p_i$ for $i\in J$, so $t_i=\pm 1$. If $\mathop{\Sigma}\limits_{i\in J}\,t_i\ge\left|J\right|-1$ then $t_i=1$ and $X_i=p_i$ for $i\in J$ and so the preceding sum is actually $\left|J\right|$. This proves A): $S_J\cap \mathcal{P}=\emptyset$. If, further, $(X,Y)\in \mathcal{P_J}$ then $X_i=p_i$ also for $i\notin J$ so $X=p$ and we have proved a part of B): $D_J\cap \mathcal{P}_J=\{p\}$.   

\bigskip

%To see that $S_J$ lies in $\mathcal{E}$ just observe that a point in the intersection $S_J \cap\mathcal{P}$ must satisfy for $i\in J$ both
%$X_i=x_i e^i_0$ and $ X_i\cdot X_i = \frac{1}{k}$, so $x_i= \frac{\pm 1}{\sqrt{k}}$. But this implies that

%$$\mathop{\Sigma}\limits_{i\in J}\,x_i=\frac{q}{\sqrt{k}}.$$

%\noindent where $q$ is congruent to $\left|J\right|$ modulo 2, which is incompatible with the last equation for $S_J$ and the intersection is empty.

%, since  From the equations of $S_J$ we obtain that $$X_i=p_i$$ \noindent also for $i\notin J$ and that $Y=0.$Therefore the only point in $D_J \cap\mathcal{P}_J$ is $$p=\mathop{\Sigma}\limits_{i=1}^k\frac{1}{\sqrt{k}}e^i_0=(p_1,\dots,p_k,0).$$

To verify that the intersection is transversal just observe that the normal space to $\mathcal{P}_J$ at $p$ is exactly $L_J$. Therefore, $\mathcal{P}_J$ and $L_J$ intersect transversely at $p$ in $\R^{n+1}$. It follows that $\mathcal{P}_J$ and $L_J\cap S^n$ also intersect transversely in $S^n$ at $p\in D_J$ and B) is proved.

%$p$ points the unit normal vector to $S^n$ at $p$ is $\mathop{\Sigma}\limits_{i=1}^k e^i_0,$ so a tangent vector $v$ to $S^n$ at $p$ is characterized by having the sum of its $e^i_0$ components equal to zero. We can decompose $v$ into a sum $v_1+v_2$ where $v_1$ is the projection of $v$ on the subspace of $\R^n$ generated by the vectors $e^i_j$ where $i\in J$ and $j>0$. Such a vector$e^i_j$ is tangent to the sphere $S^{n_i}$ at the point $p_i$, so $v_1$ is tangent to the product $\mathcal{P}_J$ at $p$. Now $v_2$ has all those components equal to zero and the sum of its $e^i_0$ components equal to zero. But this fact clearly characterizes the tangent vectors to $D_J$ at $p$. So $v_2$ is tangent to $D_J$ at $p$ and therefore the intersection of $\mathcal{P}_J$ and $D_J$ is transversal at $p$ since the (direct) sum of the two tangent spaces is the whole tangent space of $S^n$. So the linking number of $\mathcal{P}_J$ and $S_J$ is $\pm 1$.

\bigskip

To prove C) in each dimension first observe that if $J\subset J'$ they would be equal since $n_J=n_{J'}$ and by hypothesis every $n_i$ is positive. Then there is an $i\in J$ not in $J'$ so $X_i=p_i$ both in $D_J$ and in $\mathcal{P}_{J'}$. Now if we take $q_i \in \R^{n_i+1}$ close to $p_i$ with $q_i\ne p_i$ and $q_i\cdot q_i=k$ (which is possible because $n_i$ is positive) we get a new product $\mathcal{P}'_{J'}$ by changing the equation $X_i=p_i$ into $X_i=q_i$. So we have $D_J \cap\mathcal{P}'_{J'}=\emptyset$. Since $\mathcal{P}'_{J'}$ is close and isotopic to $\mathcal{P}_{J'}$ we have proved C).

\bigskip

So the linking number of $S_J$ with $\mathcal{P}_{J'}$ is zero and Theorem A2.1 is proved.

\bigskip

For $k>1$ we cannot apply to $\mathcal{E}$ Theorem A.1 characterizing a connected sum along the boundary (and in fact it is not, since its boundary is a product of at least
three spheres). But we can apply it to $\mathcal{E}\times D^1$ since the spheres $S_J$ can be embedded with different values of the coordinate in $D^1$ so they are disjoint,
represent a basis of its homology and have codimension at least 3 and dimension at least 2 if $q+k\ge 4$. As a consequence we have:

\bigskip

\textbf{Theorem A2.2.} \textit{If $q+k\ge 4$ , then}

\bigskip

A) $\mathcal{E}\times D^1$ \textit{is diffeomorphic to} $$\scf_{\emptyset\ne J \subset I} (S^{n-n_J-1}\times D^{n_J+2}).$$

B) The double of $\mathcal{E}$, $\mathcal{D} (\mathcal{E})$, \textit{is diffeomorphic to} $$\sc_{\emptyset\ne J \subset I} (S^{n-n_J-1}\times S^{n_J+1}).$$

C) The trivial open book on $\mathcal{E}$, $\mathcal{TOB} (\mathcal{E})$, \textit{is diffeomorphic to} $$\sc_{\emptyset\ne J \subset I} (S^{n-n_J-1}\times S^{n_J+2}).$$

\bigskip

If $n_i>0$ for $i=1,\dots , k$ then A) follows from Theorem A1 since the hypothesis $q+k\ge 4$ implies that $\mathcal{P}$ has codimension at least 3 in $S^n$ and so
$\mathcal{E}\times D^1$ is simply connected with simply connected boundary. When some $n_i$ are zero it follows by induction on the number of them: if we add a new
factor $S^0$ to a product $\mathcal{P}$ we get a new complement $\mathcal{E}'$ which is clearly the connected sum of two copies of $\mathcal{E}$. By Lemma 1.2, $\mathcal{E}'\times D^1$
is the connected sum along the boundary of two copies of $\mathcal{E}\times D^1$ and one copy of $S^{n-1} \times D^2$ and by the induction hypothesis the collection of summands
of the form $S^{n-n_J-1}\times D^{n_J+2}$ coincides with the one claimed for $\mathcal{E}'\times D^1$ and A) is proved. B) follows from A) by taking its boundary and C) by taking its
product with $D^1$ and again the boundary.

\bigskip

A version of Part B) was proved in the PL-category by McGavran ([Mc]) for the case where all $n_i=1$ (i.e., when $\mathcal{P}$ is the $k$-torus) with a possibly more general hypothesis on the embedding. It does not seem that analogs of the other parts of the theorem can be deduced easily from his results.

For our purpose (the proof of Theorem 2.2) we need a different version in the torus case, where the standard nature of the embedding is obtained from a codimension condition:

\bigskip

\textbf{Theorem A2.3.} \textit{Consider an embedding of the torus $T^k$ in $S^n$ with $2k<n$ and $\mathcal{E}$ the complement of an open tubular neighborhood of $T^k.$ Then}

A) $\mathcal{E}\times D^1$ \textit{is diffeomorphic to} $$\scf_{j=1}^k {k\choose j} (S^{n-j-1}\times D^{j+2}).$$

B) \textit{The double of} $\mathcal{E}$, $\mathcal{D} (\mathcal{E})$, \textit{is diffeomorphic to} $$\sc_{j=1}^{k} {k\choose j} (S^{n-j-1}\times S^{j+1}).$$

C) \textit{The trivial open book on} $\mathcal{E}$, $\mathcal{TOB} (\mathcal{E})$, \textit{is diffeomorphic to} $$\sc_{j=1}^{k}
{k\choose j} (S^{n-j-1}\times S^{j+2}).$$

\bigskip

This follows from the previous theorem since any embedding of the torus $T^k$ in $S^n$ with $2k<n$ is isotopic to the product considered there.

\bigskip

Obviously a similar result holds for any embedding of the above product of spheres if $q+k-1>\mathop{\Sigma}\limits_{i=1}^k n_i$, with refinements depending on the minimum of the the $n_i$.

\newpage

\section*{References}

[A-LL] D. Allen, and J. La Luz, \emph{A Counterexample to a conjecture of
Bosio and Meersseman}, Contemporary Mathematics 460, AMS, 2008, 37-46.

\bigskip

\noindent[B-B-C-G], A. Bahri, M. Bendersky, F. R. Cohen, and S. Gitler, \emph{The polyhedral product functor: a
method of computation for moment-angle complexes, arrangements and related spaces}, Advances in Mathematics
Volume 225 (2010), 1634-1668.

\bigskip

\noindent[B-B-C-G2], A. Bahri, M. Bendersky, F. R. Cohen, and S. Gitler, \emph{Operations on polyhedral products and a new topological construction of infinite families of toric manifolds}, arXiv:1011.0094. This is a very enriched version of the 2008 draft \emph{An infinite family of toric manifolds associated to a given one}.

\bigskip

\noindent[B-B-C-G3], A. Bahri, M. Bendersky, F. R. Cohen, and S. Gitler, \emph{Cup-products for the polyhedral product functor}, to appear in Math. Proc. Camb. Phil Soc.

\bigskip

\noindent[Ba1] I.V. Baskakov, \emph{Cohomology of $K$-powers of spaces and the combinatorics of simplicial subdivisions}
Russian Math. Surveys 57(2002) 989-990.

\bigskip

\noindent[Ba2] I.V. Baskakov, \emph{Massey triple products in the cohomology of moment-angle complexes.} Russian Math. Surveys 58 (2003), no.5,
1039-1041.

\bigskip

\noindent[B], F. Bosio, \emph{Vari\'et\'es complexes compactes: une g\'en\'eralisation de la \newline construction de Meersseman
et L\'opez de Medrano-Verjovsky}. Ann. Inst. Fourier (Grenoble), 51 (2001), 1259-1297.

\bigskip

\noindent[B-M], F. Bosio and L. Meersseman, \emph{Real quadrics in $\C^n$, complex manifolds and convex polytopes}. Acta Math. 197 (2006), no. 1, 53-127.

\bigskip

\noindent[Br], W. Browder, \emph{Surgery on simply-connected manifolds.} Ergebnisse der Mathematik und ihrer Grenzgebiete, Band 65. Springer-Verlag, New York-Heidelberg, 1972. ix+132p.

\bigskip

\noindent[B-P], V.M. Buchstaber and T.E. Panov, \emph{Torus actions and their applications in Topology and
Combinatorics}, University Lecture Seriers, AMS (2002).

\bigskip

\noindent[C-K-P], C. Camacho, N. Kuiper and J. Palls. \emph{The topology of holomorphic
flows with singularities}, Publications Mathematiques I.H.E.S., 48 (1978), 5-38.

\bigskip

\noindent[Ch], M. Chaperon, \emph{G\'eom\'etrie diff\'erentielle et singularit\'es de syst\`emes dynamiques.} Ast\'erisque 138-139 (1986),
434 pages.

\bigskip

\noindent[Ch-LdM], M. Chaperon and S. L\'opez de Medrano,
\emph{Birth of attracting \newline compact invariant submanifolds
diffeomorphic to moment-angle manifolds in generic families of
dynamics}, C. R. Acad. Sci. Paris, Ser. I 346 (2008) 1099-1102.

\bigskip

\noindent[dL], M. de Longueville, \emph{The ring structure on the cohomology of coordinate subspace arrangements},
Math. Z. 233 (2000) 553-577.

\bigskip

\noindent [D-S], G. Denham, A. Suciu, \emph{Moment-angle complexes, monomial ideals, and Massey products}, Pure and
Applied Math. Quarterly 3 (2007), no. 1, 25-60.

\bigskip

\noindent[D-J], M. Davis and T. Januszkiewicz, \emph{Convex polytopes, Coxeter orbifolds and torus actions}, Duke
Math. Journal 62 (1991), 417-451.

\bigskip

\noindent[G-P-W], V. Gasharov, I. Peeva and V. Welker, \textit{Coordinate Subspace Arrangements and Monomial Ideals},
Adv. Stud. Pure Math., 33, Math. Soc. Japan, Tokyo, (2002) 65–74. 

\bigskip

\noindent[GA], F. Gonz\'alez Acu\~na, \emph{Open Books}, University of Iowa Notes, 12 p.

\bigskip

\noindent[G-LdM], S. Gitler and S. L\'opez de Medrano, \emph{Cup products in $(D^1,S^0)$ generalized moment angle complexes}, in preparation.

\bigskip

\noindent[H], F. Hirzebruch, private conversation, based on \emph{Arrangements of lines and algebraic surfaces}, in
Arithmetic and Geometry, Vol. II (= Progr. Math. 36), Boston: Birkhauser, 1983, 113140.

\bigskip

\noindent[K-M] M. Kervaire and J. Milnor, \emph{Groups of Homotopy Spheres: I}
Annals of Mathematics, 2nd Ser., Vol. 77, No. 3. (1963), pp. 504-537.

\bigskip

\noindent[L-M] J. J. Loeb and M. Nicolau, \emph{On the complex
geometry of a class of non-Kaehlerian manifolds.} Israel J. Math.
110  (1999), 371--379.

\bigskip

\noindent[LdM1], S. L\'opez de Medrano, \emph{Topology of the intersection of quadrics in $\R^n$}, in {\it Algebraic Topology}
(Arcata Ca, 1986), Springer Verlag LNM  \textbf{1370}(1989), pp. 280-292, Springer Verlag.

\bigskip

\noindent[LdM2], S. L\'opez de Medrano, \emph{The space of Siegel leaves of a holomorphic vector field}, in
Holomorphic Dynamics (Mexico, 1986), Lecture Notes in Math., 1345, pp. 233-245. Springer, Berlin, 1988.

\bigskip

\noindent[LdM3], S. L\'opez de Medrano, \emph{Intersections of quadrics and the polyhedral product functor}, in preparation.

\bigskip

\noindent[LdM-V], S. L\'opez de Medrano and A. Verjovsky, \emph{A new family of complex, compact, nonsymplectic
manifolds.} Bol. Soc. Brasil. Mat., 28 (1997), 253-269.

\bigskip

\noindent[Mc] D. McGavran, \emph{Adjacent connected sums and torus actions},
Trans. AMS 251 (1979) 235-254

\bigskip

\noindent[Me], L. Meersseman, \emph{A new geometric construction of compact complex manifolds in any
dimension}. Math. Ann., 317 (2000), 79-115.

\bigskip

\noindent[Me-V1], L. Meersseman and A. Verjovsky, \emph{Holomorphic principal bundles over projective toric
varieties.} J. Reine Angew. Math., 572 (2004), 57-96.

\bigskip

\noindent[Me-V2], L. Meersseman and A. Verjovsky, \emph{Sur les vari\'et\'es LV-M}, Contemporary Mathematics 475, AMS, 2008, 111-134.

\bigskip

\noindent[W1], C.T.C. Wall, \emph{Classification of $(n-1)$-connected $2n$-manifolds.}  Ann. of Math. (2)  75  1962 163--189.

\bigskip

\noindent[W2], C.T.C. Wall, \emph{Classification problems in differential topology. VI. Classification of $(s-1)$-connected $(2s+1)$-manifolds.}  Topology  6  1967 273--296.

\bigskip

\noindent[W3], C.T.C. Wall, \emph{Stability, pencils and polytopes}, Bull. London
Math. Soc., 12 (1980), 401-421.

\bigskip

samuel.gitler@gmail.com

Department of Mathematics, Cinvestav, Mexico.

El Colegio Nacional.

\bigskip

santiago@matem.unam.mx

Instituto de Matem\'aticas, Universidad Nacional Aut\'onoma de M\'exico.

\end{document}